\documentclass[11pt]{article}
\usepackage{amsthm}
\usepackage{amsmath,amsfonts,amssymb}
\usepackage{graphicx, color}

\def\bega{\begin{array}}
\def\enda{\end{array}}
\def\ds{\displaystyle}
\def\begi{\begin{itemize}}
\def\endi{\end{itemize}}

\def\argmin{\hbox{arg}\!\min}
\def\ve{\varepsilon}

\def\J{{\cal J}}

\def\L{{\bf L}}

\def\R{{I\!\!R}}
\def\bfv{{\bf v}}
\def\bfw{{\bf w}}
\def\bpm{\begin{pmatrix}}
\def\epm{\end{pmatrix}}

\def\implies{\Longrightarrow}
\def\vp{\varphi}

\def\vs{\vskip 2em}

\def\v{\vskip 1em}
\def\wto{\rightharpoonup}
\def\C{{\cal C}}

\def\D{{\cal D}}

\def\ov{\overline}
\def\forall{\hbox{for all}~}
\def\Tilde{\widetilde}
\def\Hat{\widehat}
\def\meas{\hbox{meas}}
\def\bel{\begin{equation}\label}
\def\eeq{\end{equation}}
\def\sqr#1#2{\vbox{\hrule height .#2pt
\hbox{\vrule width .#2pt height #1pt \kern #1pt
\vrule width .#2pt}\hrule height .#2pt }}
\def\square{\sqr74}
\def\endproof{\hphantom{MM}\hfill\llap{$\square$}\goodbreak}

\newtheorem{theorem}{Theorem}[section]

\newtheorem{lemma}{Lemma}[section]

\newtheorem{definition}{Definition}[section]
\newtheorem{example}{Example}[section]

\begin{document}

\title{\bf Controlled Traveling Profiles for Models of Invasive Biological Species}
\vs

\author{Alberto Bressan and Minyan Zhang\\
\,
\\
Department of Mathematics, Penn State University \\
University Park, Pa.~16802, USA.\\
\,
\\
e-mails: axb62@psu.edu, muz94@psu.edu
}
\maketitle

\begin{abstract} We consider a family of controlled reaction-diffusion equations, 
describing the spatial spreading of an invasive biological species.  
For a given propagation speed $c\in\R$,
we seek a control with minimum cost, which achieves a traveling profile with speed $c$.
For various nonlinear models, the existence of a (possibly measure valued) 
optimal control is proved, together with necessary conditions for optimality.
In the last section we study a case where the wave speed cannot be modified by 
any control with finite cost.

The present analysis is motivated by the recent results in \cite{BCS1, BCS2}, 
showing how
a control problem for a reaction-diffusion equation can be approximated
by a simpler problem of optimal control of a moving set.
\end{abstract}

\section{Introduction}
\label{sec:1}
\setcounter{equation}{0}

Consider a reaction-diffusion equation
of the form
\bel{1} u_t~=~\sigma\Delta u + f(u,\alpha).\eeq
Here $t\geq 0$ is time, $x\in \R^n$ is the spatial variable, while $u=u(t,x)$ denotes the density
of an invasive biological species, such as mosquitoes.  
We assume that, by implementing a control $\alpha=\alpha(t,x)\geq 0$,
the population can be partly removed. This will slow down, or even reverse, its spatial propagation.

By a rescaling of the dependent variable we shall always assume that, when $\alpha=0$, i.e.~in absence of control, 
one has
$$f(0,0)~=~f(1,0)~=~0.$$
In other words, the maximum population density sustained by  the environment (i.e., the carrying capacity)   is normalized so that $u^{max}=1$. 

Given an initial density 
\bel{ipd}u(0,x)\,=\, \bar u(x)\eeq
and a time interval $[0,T]$, 
a natural objective can be stated as
\bel{min1}\hbox{minimize:}\quad 
\J~\doteq~\int_0^T \left( \int_{R^n} \bigl[u(t,x) +\alpha(t,x) \bigr]dx\right)dt\,.
\eeq
%for suitable initial and boundary conditions.
The right hand side of (\ref{min1})
accounts for the  population size, plus the cost of the control,  integrated over time.
%
%A second model we shall study  includes
%the loss of trees, which get contaminated by a pathogen carried by the insects.
%Calling $\theta=\theta(t,x)\in [0,1]$ the fraction of trees that are contaminated,
%this leads to
%\bel{tht}
%\theta_t~=~ \kappa_2 u (1-\theta).\eeq
%Here  the constant $\kappa_2$ is an infection rate.
%The corresponding optimization problem takes the form
%\bel{min2}  \hbox{minimize:}\quad 
%\J~=~\int_0^T \left( \int_\Omega \bigl[ \alpha(t,x)+\theta(t,x) \bigr]dx\right)dt. \eeq
%\v
%%While the existence of an optimal control strategy can be inferred from general principles,
%%it is often very difficult to explicitly describe the optimal solutions. 

Thanks to the fact that $u=0$ and $u=1$ are equilibrium states,  in many cases the  solution  to (\ref{1})
can be approximately 
described in terms of the set $\Omega(t)$ where $u(t,x)\approx 1$.
 Namely, if the diffusion coefficient $\sigma>0$  is small,  we expect that the
difference $\bigl\| \chi_{\Omega(t)} - u(t,\cdot)\bigr\|_{\L^1}$ will also be small. The
 characteristic function $\chi_{\Omega(t)} $  of the set $\Omega(t)$ thus
provides a good approximation to the density $u(t,\cdot)$ itself. 
 Based on this observation, in \cite{BCS1} it was proposed to replace the  problem 
(\ref{min1}) by an optimization problem for the moving set $\Omega(t)$.
More precisely, let $c(t,x)$ be the speed at which the boundary of the set $\Omega(t)$ moves, in the direction of the interior normal, at a point $x\in \partial\Omega(t)$.
The new optimization problem then takes the form
\bel{minom}
\hbox{minimize:}\quad 
\J~=~\int_0^T \left( \meas \bigl(\Omega(t)\bigr)+ \int_{\partial\Omega(t) } E\bigl(c(t,x)\bigr)\, dx \right) dt.
\eeq
The cost function $E(c)$, which is integrated over the boundary of the set $\Omega(t)$, 
measures the effort  needed to push the boundary inward with speed $c$. 
As shown in \cite{BCS1}, it is this particular function that provides the link between the two problems (\ref{minom}) and (\ref{min1}). A rigorous justification of this approximation procedure can be achieved via a sharp interface limit.

With this motivation in mind, in the present paper we study the function 
$E(\cdot)$, for various nonlinear 
parabolic equations, or systems.
In our basic setting, $E(c)$ is defined as the minimum cost of a control $\alpha(\cdot)$ which yields a traveling wave solution to 
(\ref{1}) with speed $c$.  This leads to the problem
\bel{mina}
\hbox{minimize:}\quad \|\alpha\|_{\L^1}\eeq
among all integrable functions $\alpha\geq 0$ such that there exists a solution to the 
ODE
\bel{ode0}\sigma
U'' + c \,U' +f(U,\alpha)~=~0,\eeq
with asymptotic conditions
\bel{ac0}
U(-\infty)~=~0,\qquad\qquad U(+\infty)~=~1.\eeq
%Denoting by $E(c)$ the minimum cost for the problem (\ref{mina})--(\ref{ac0}),
%the analysis in  \cite{BCS1} has established a basic connection between the optimization problems
%(\ref{min1}) and (\ref{minom}).   
%
In the models considered in \cite{BCS1}, the function $f(u,\alpha)$ has linear dependence on 
the control variable $\alpha$.  
Namely, the two main cases
$$f(u,\alpha)~=~F(u) + \alpha,\qquad\hbox{or}\qquad f(u,\alpha)~=~F(u)+\alpha u,$$
were studied.
Thanks to this assumption, the difference in cost between any two admissible controls
can be directly computed by Stokes' formula \cite{BoP, HH}. This yields a straightforward 
way to identify the optimal solution.

In the present paper, our first goal is to extend the analysis of controlled traveling waves
to a more general class of functions $f$, possibly nonlinear also w.r.t.~the control variable $\alpha$.
In this case the techniques from \cite{HH} cannot be implemented, 
and the construction of optimal profiles requires a more careful analysis.

In the second part of the paper,
 we focus our study on two systems of PDEs, describing the interaction 
between disease-carrying insects and infected trees. 
A relevant example is provided by
{\it Xylella fastidiosa}, which is a  plant pathogenic bacterium that attacks olive trees.
It is transmitted by a
meadow spittlebug, the {\it Philaenus spumarius}, a sap-feeding insect.
In \cite{ACS} a detailed model for spatial propagation of a Xylella was introduced.
This is described by a system of four equations for the densities of
(i) healthy and 
infected insects, and (ii)
healthy and infected trees.
Here we consider two simplified models, that will allow a more detailed mathematical analysis.
\v
{\bf Model 1.} Assume that:
\begi
\item The insect population spreads by diffusion and reproductive growth.
\item By spraying pesticides, some of the insects can be removed.  
This slows down, or even reverses, their spatial propagation.
\item
All insects carry the infection, and contaminate the trees.
\endi
Calling
\begi\item $u=u(t,x)\in [0,1]$ the density of insects,
\item $\theta=\theta(t,x)\in [0,1]$ the fraction of trees that are infected,
\item $\alpha=\alpha(t,x)\geq 0$ the control function,
\endi
the evolution of these variables can be described by
\bel{ode1}
\left\{ \bega{rl} u_t&=~\Delta u + f(u,\alpha)\,,\\[2mm]
\theta_t&=~ \kappa_1 u (1-\theta)\,.\enda\right.\eeq
Here  the constant $\kappa_1$ is an infection rate.
The function $f=f(u,\alpha)$, modeling the controlled population growth, can
take different forms.  For example:
\begi
\item[(i)] Logistic growth + insect removal by pesticides or mosquito nets. This leads to
\bel{f3} 
f(u,\alpha) ~=~\kappa_3 (1-u) u  - \alpha u.\eeq
\item[(ii)] Weed removal, reducing the carrying capacity of the ecosystem.
A possible model is
\bel{f4} 
f(u,\alpha) ~=~ u(u-u^*)\left[ {1\over 1+\alpha} - u\right],\eeq
where $u^*\in [0, 1/2]$.
Notice that in this case the maximum population supported by the environment
shrinks to $(1+\alpha)^{-1}<1$
as the control $\alpha$ increases.  This is another way to reduce the density of insects.
\endi

For the above model, a natural goal is to minimize
\bel{min13} \bega{rl}
\J&\ds=~\int_0^T \left( \int_{\R^2} \bigl[ \alpha(t,x)+\theta(t,x) \bigr]dx\right)dt\\[4mm]
&=~\hbox{[cost of the control]} + \hbox{[fruit production loss]}.
\enda
\eeq
for given initial data.

\v
{\bf Model 2.} We here assume that 
\begi
%\item The insect population is already spread out over the whole space. 
%However, only a subregion is occupied by infected insects.
\item Newly born insects are healthy. Only later in life they can be infected,
by the presence of contaminated trees. 
\item Infected insects contaminate the trees, and contaminated trees infect the new insects.
\item By spraying pesticides, some of the insects can be removed.  
%By lowering the density of  insects near the interface between healthy and infected ones, 
%we  create a wedge between the two populations, thus reducing the speed at which the 
%contamination spreads in the territory.
\endi
In addition to the previous variables, calling
\begi
%\item $u=u(t,x)\,\geq\,0$ the density of insects,
\item $I=I(t,x)\in [0,1]$ the fraction of insects which are infected,
\item $v=Iu$ the density of infected insects,
%\item $\theta=\theta(t,x)\in [0,1]$ the fraction of trees that are contaminated,
%\item $\alpha=\alpha(t,x)\geq 0$ a control function.
\endi
we thus consider the system of evolution equations
%For simplicity, we restrict the analysis to a case where the control enters linearly in the equations.  
%The evolution of these variables is now governed by
\bel{ode2}
\left\{ \bega{rl} u_t&=~\Delta u + f(u) - \alpha u\,,\\[2mm]
(Iu)_t &=~ \Delta(Iu) + \kappa_2 (1-I)u  \theta  - \alpha I u - d \,I u,\\[2mm]
\theta_t&=~ \kappa_1  Iu (1-\theta) \,.\enda\right.\eeq

The constants $\kappa_1, \kappa_2$ are infection rates, while $d$ is a death rate.   

Motivated by \cite{BCS1}, for the three models (\ref{1}), (\ref{ode1}), (\ref{ode2}), 
we are interested in (i) the existence of controlled traveling profiles having a given speed $c$, and
(ii) control functions $\alpha(\cdot)$ which achieve these traveling profiles and have minimum cost.

\begin{figure}[ht]
\centerline{\hbox{\includegraphics[width=11cm]{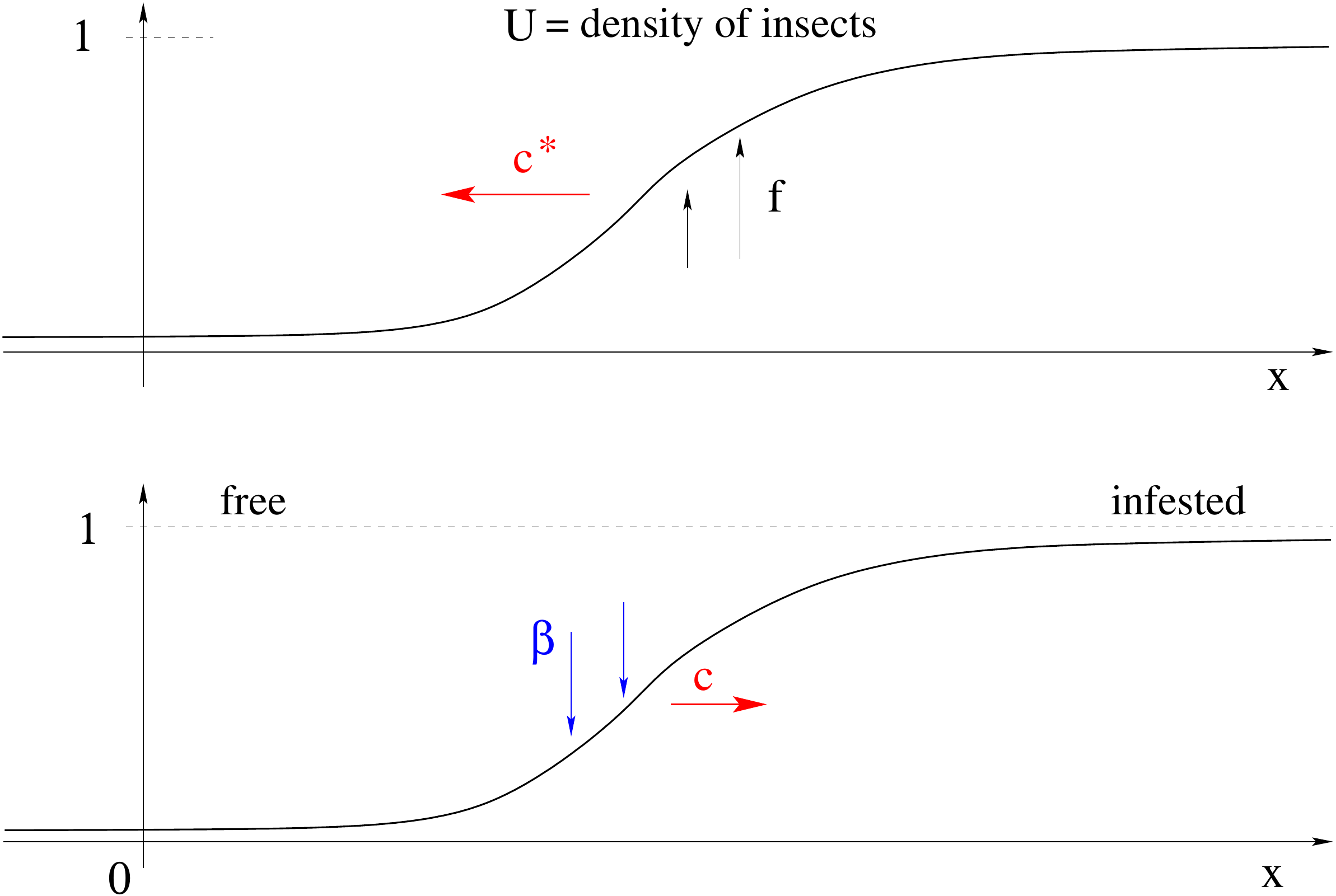}}}
\caption{\small Traveling profiles for Model 1.  Above: without any control, the insect population spreads toward the left, with a speed $c^*<0$.
Below: applying a control, part of the population is removed.   This yields a new traveling wave profile,
with speed $c>c^*$.  % When $c>0$, as time progresses the population is reduced. 
}
\label{f:tw1}
\end{figure}

We now summarize the main results, proved in the remainder of the paper.
In Section~\ref{sec:2} we study the scalar equation (\ref{1}). 
By a rescaling of the spatial variable, it is
not restrictive to assume $\sigma=1$.  In absence of control,
by the standard theory in \cite{F, VVV}
it is known that the 
equation admits a traveling wave solution with a suitable speed $c^*<0$.  
Here we prove that, given any speed $c>c^*$,
there exists a control function $\alpha(\cdot)$ with finite cost which yields a traveling profile with speed $c$. 
More precisely (see Fig.~\ref{f:tw1}), setting
$$u= U(x-ct),\qquad \alpha = \alpha(x-ct),$$
we construct a solution to 
\bel{TW0}U'' + c U' +f(U,\alpha)~= 0,\eeq
with asymptotic conditions (\ref{ac0}).
%\bel{UT0}
%U(-\infty)\,=\,0,\qquad\qquad U(+\infty) = 1.
%\eeq
 In Section~\ref{sec:3} we prove that a suitable control
function $\alpha^*(\cdot)$ can be chosen, having minimum cost.  
Necessary conditions for optimality are then 
derived in Section~\ref{sec:4}. 
In turn, these can be used in a shooting method, to numerically compute optimal solutions.   
Plots of an optimal  traveling profile, and of the
minimum cost  $E(c)$ as a function of the speed $c$, are shown in Fig.~\ref{f:c01profile} and 
Fig.~\ref{f:minco}, respectively.

In Section~\ref{sec:5} we study Model~1. Here the main result shows that,
for every wave speed $c\in [c^*, 0[$, the system (\ref{ode1}) admits a controlled traveling wave
with speed $c$. In other words, by removing part of the pest population, the speed at which the 
contamination advances can be slowed down to almost  zero.

The last two sections are concerned with Model 2.
Looking for traveling wave solutions of (\ref{ode2}) of the form
$$u(t,x)~=~U(x-ct),\qquad I(t,x)~=~I(x-ct),\qquad \theta(t,x)~=~\Theta(x-ct),$$
we are led to the system of three ODEs:
\bel{TW2} 
\left\{ \bega{rl} U'' + c U' +f(U)-  \alpha U&=~ 0,\\[2mm]
(IU)'' + c (I'U+IU')  + \kappa_2 (1-I) U\Theta-\alpha IU - d\, IU&=~0 ,\\[2mm]
c\Theta' + \kappa_1 (1-\Theta) I U  &= ~0.\enda\right.\eeq
Two scenarios can be considered.   In Section~\ref{sec:6}
we study (\ref{TW2}) with asymptotic conditions
\bel{UT2}\left\{ \bega{rl}U(-\infty)&=~ 0,\\[1mm]
I(-\infty)&=~ 0,\\[1mm]\Theta(-\infty)&=~ 0,
\enda\right.\qquad   \qquad 
\left\{ \bega{rl}U(+\infty)&=~ 1,\\[1mm]
I(+\infty)&=~ I^*,\\[1mm]\Theta(+\infty)&=~ 1,
\enda\right.\eeq
where $I^*= \kappa_2/(\kappa_2+d)$.
In other words, the density of insects is  vanishingly small as $x\to -\infty$, but large for $x\to +\infty$.
All trees are healthy in the limit as $x\to -\infty$, while they are increasingly infected as $x \to +\infty$.   In this case, controlling the contamination
essentially amounts to slowing down the spreading of the insect population (see Fig.~\ref{f:tw1}). 
Observing that the density of infected insects trivially satisfies $IU\leq U$, by 
a comparison argument we prove that, if 
the control $\alpha=\alpha(x-ct)$ yields a traveling profile with speed
$c<0$  for the first equation in (\ref{TW2}), then 
the same control yields a traveling profile for the entire system (\ref{TW2}), with the same speed.

\begin{figure}[ht]
\centerline{\hbox{\includegraphics[width=11cm]{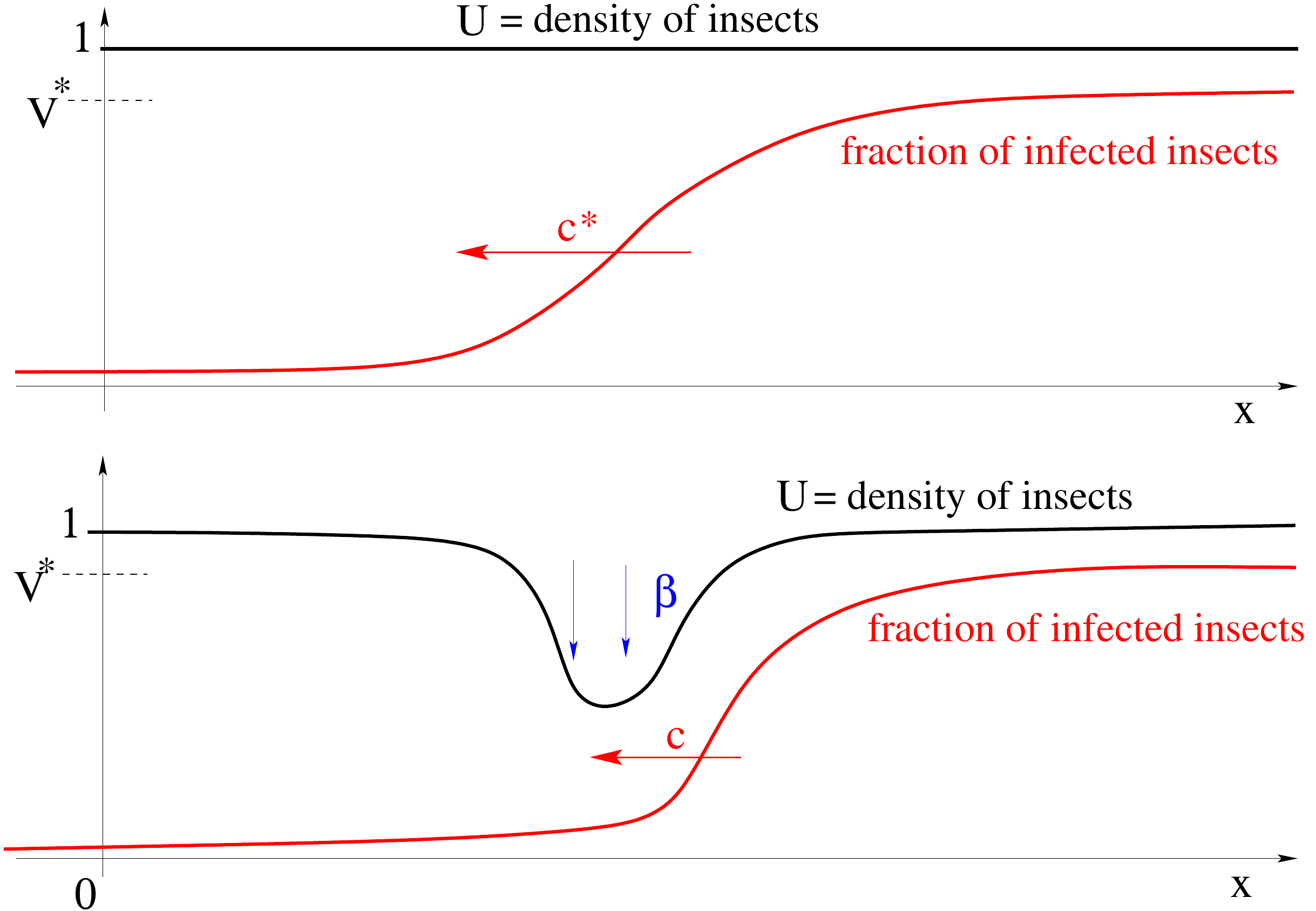}}}
\caption{\small Traveling profiles for Model 2.  Above:  without any control, the insect population reaches everywhere its
maximum value $U=1$, while the fraction of infected insects keeps increasing, propagating to the left with speed $c^*<0$.
 Below: applying a control, part of the population is removed, in a neighborhood of the interface between
 healthy and infected individuals.   This yields a different traveling wave profile.
 However, our analysis shows that the propagation speed cannot be affected. }
\label{f:tw2}
\end{figure}

Finally, in Section~\ref{sec:7} we consider again the system (\ref{TW2}), but with asymptotic 
conditions
\bel{UT3}\left\{ \bega{rl}U(-\infty)&=~ 1,\\[1mm]
I(-\infty)&=~ 0,\\[1mm]\Theta(-\infty)&=~ 0,
\enda\right.\qquad   \qquad 
\left\{ \bega{rl}U(+\infty)&=~ 1,\\[1mm]
I(+\infty)&=~ I^*,\\[1mm]\Theta(+\infty)&=~ 1.
\enda\right.\eeq
Notice that here
the density of insects is large  for $x\to +\infty$  as well as for $x\to -\infty$. 
Insects and trees
are all healthy in the limit as $x\to -\infty$, while they are increasingly infected as $x \to +\infty$.

In the uncontrolled case where $\alpha=0$, one would have a traveling wave profile where the insect population is everywhere constant:
$U(x)=1$. On the other hand, as shown at the top of Fig.~\ref{f:tw2},  the fraction of infected  trees and insects keeps increasing. Indeed, the contamination advances toward the left, with speed $c^*<0$.

An interesting question now arises.   Assume that, by applying a control, we locally 
reduce the population
density $U$.   As shown at the bottom of Fig.~\ref{f:tw2}, this will create a buffer between a region (to the right)
where most of the trees and insects are infected, and a region (to the left)
where trees and insects are still largely healthy.  Can this strategy effectively reduce the speed
at which the contamination advances\,?

Our analysis shows that the answer is negative.  Indeed, the speed of a traveling wave must
satisfy a constraint stemming from the linearization of the system (\ref{TW2}) at
the asymptotic state
$(U,I,\Theta)=(1,0,0)$.  We now observe that any control $\alpha(\cdot)$ with finite cost must be integrable, hence vanishingly small as $x\to -\infty$.
As a consequence, the presence of this additional control cannot
remove the above constraint on the wave speed.   A precise statement of the result is given in Theorem~\ref{t:71}. 
\v
Traveling profiles for systems of parabolic equations is a classical subject, with an extensive literature. 
See for example \cite{F,  HI, LZ, VVV} and references therein.
Control problems for nonlinear parabolic equations, such as optimal harvesting problems,
were studied in \cite{CG, CGS, LM, RBZ}.  For more accurate models of the spreading and control
of invasive populations we refer to
\cite{ACD, ACM, ACS, CRo, SBGW}.   
%
%According to the analysis in \cite{BCS1}, by determining
%the minimum cost  $E(c)$ of a traveling wave with speed $c$, one can then approximate the original 
%problem (\ref{min1}) with an optimization problem of the form (\ref{minom}) for a moving set.
Given an effort function $E(c)$, optimization problems  for a moving set 
of the form (\ref{minom}) have been recently studied in \cite{BCS2},
proving the existence of optimal strategies and establishing necessary conditions for optimality.
Control problems for a moving set, describing the support of a population, have also been
considered in
\cite{BZ, CPo, CLP}.

\section{Controlling a traveling front}
\label{sec:2}
\setcounter{equation}{0}
Given $c\in\R$, as in (\ref{mina})--(\ref{ac0}) we seek a control $\alpha(\cdot)$ with minimum $\L^1$ norm, that 
produces a traveling wave with speed $c$.
Assuming for simplicity that $\sigma=1$, and using the notation
\bel{ba}\beta\,=\, f(u,0)- f(u,\alpha),\qquad\qquad f(u)\,=\,f(u,0),\eeq
we can write 
(\ref{1}) in the form
\bel{fu1} u_t ~=~\Delta u + f(u) -\beta.\eeq
In addition,  we introduce the cost function $L$ implicitly defined by
\bel{Lag}
 L(u,\beta)\,=\,\alpha.\eeq
The optimization problem for traveling wave profiles can now be stated as follows.
\begi
\item[{\bf (OTW)}] {\it Given functions $f(u)$ and $L(u,\beta)$, and a speed $c\in\R$,
find a nondecreasing profile $U:\R\mapsto [0,1]$ and a  control function 
$\beta:\R\mapsto\R_+$ which minimize the cost
\bel{cp1}
J(U,\beta)~\doteq~ \int_{-\infty}^{+\infty} L\bigl(U(x), \beta(x)\bigr)\, dx\,,\eeq
subject to
\bel{cp2}
U'' + c U' + f(U) -\beta~=~0,\qquad\qquad  U(-\infty) =0,\quad U(+\infty) =1.\eeq
}
\endi

\v

\begin{example} {\rm When $f(u,\alpha)$ is the function in  (\ref{f3}), with the notation 
introduced at (\ref{ba}), (\ref{Lag}) 
we obtain
\bel{fu2}f(u)~=~ \kappa_3(1-u) u,\qquad \beta\, =\, \alpha u,\qquad L(u,\beta)
\, =\, \alpha\, =\, {\beta\over u}\,.\eeq

On the other hand, when $f(u,\alpha)$ is the nonlinear function in  (\ref{f4}), one obtains
\bel{fu3}f(u)~=~ u (u-u^*)(1-u),\qquad \beta~=~\left( 1- {1\over 1+\alpha}\right) u(u-u^*).\eeq
%\bel{fu4}L(u,\beta)\,=\, \alpha\,=\,{\beta\over (u-u^*)u-\beta}\,.\eeq
Notice that in this case the control $\alpha$ will be effective only in the region where
$u\in [u^*,1]$, because for $u<u^*$ this control will actually increase the population growth.
As Lagrangian function, one should take
\bel{Lex}L(u,\beta)~=~\left\{ \bega{cl}  0&\qquad \hbox{if}~~\beta=0,\\[3mm]
\ds{\beta\over (u-u^*)u-\beta}&\qquad \hbox{if}~~0\leq\beta< (u-u^*) u\,,
\\[4mm]
+\infty &\qquad \hbox{in all other cases.}\enda\right.\eeq
}
\end{example}

\v
The optimization problem {\bf (OTW)}  will be studied under the following assumptions
on the  source function $f$ and the cost function $L$.
\begi
\item[{\bf (A1)}] {\it $f\in \C^2$, and moreover
\bel{f2}
f(0)~=~f(1)~=~0,\qquad f'(0)\,<\,0,\qquad f'(1)\,<\,0.\eeq
In addition, $f$ vanishes at only one intermediate point $u^*\in \,]0,1[\,$, where 
$f'(u^*)>0$.}
\item[{\bf (A2)}] {\it For every $u\in ]0,1[$ the map 
$\beta\mapsto L(u,\beta)\in \R_+\cup\{+\infty\}$ is strictly convex and has superlinear growth. More precisely, there exist constants
$C_1>0$ and $p>1$ such that
\bel{LL} L(u,0)\,=\,0,\qquad\quad L(u,\beta)~\geq~C_1\,\beta^p \qquad\quad\forall \beta\geq 0 ~\hbox{and}~u\in [0,1].\eeq
}
\endi

\begin{figure}[ht]
\centerline{\hbox{\includegraphics[width=15cm]{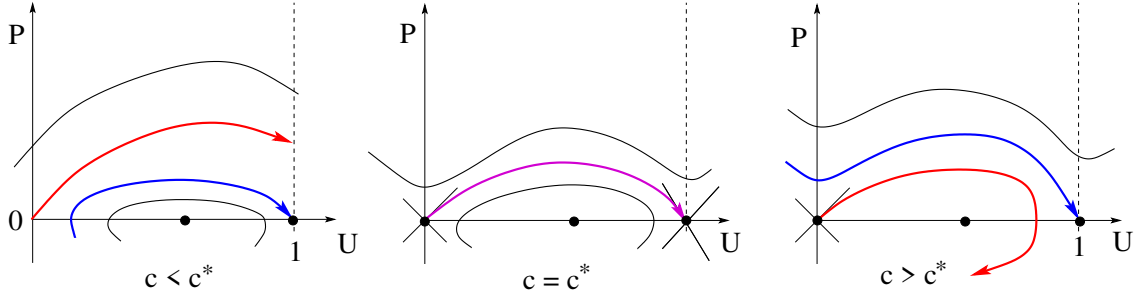}}}
\caption{\small A traveling profile for (\ref{A1}) corresponds to a heteroclinic orbit for the system
(\ref{T2}), connecting the points $(0,0)$ and $(1,0)$.  Under the assumptions {\bf (A2)}, 
such an orbit exists for one specific value $c=c^*$.}
\label{f:co5}
\end{figure}

As a preliminary, we  review some basic facts on traveling waves for reaction-diffusion 
equations of the form
\bel{A1} u_t~=~f(u) + u_{xx}\,.\eeq
By definition, a traveling profile for (\ref{A1}) with speed $c$ is a solution of the form
\bel{TP} u(t,x)~=~U(x-ct).\eeq
This can be found by solving
\bel{TE} U''+cU' + f(U)~=~0.\eeq
Assuming that $f(0)=f(1)=0$,
we seek a solution $U:\R\mapsto [0,1]$ of (\ref{TE}) with asymptotic conditions (\ref{ac0}).
%\bel{AC}
%U(-\infty) ~=~0,\qquad U(+\infty) ~=~1.\eeq
%Following \cite{F}, we set 
Setting $P=U'$,  we thus need to find a heteroclinic orbit of the system
\bel{T2}
\left\{\bega{rl} U'&=~P,\\[1mm]
P'&=~-cP-f(U),\enda\right.\eeq
connecting the equilibrium points $(0,0)$ with $(0,1)$.
A phase plane analysis of the system (\ref{T2}) yields
\begin{theorem}\label{t:21} Consider the problem (\ref{TE}), (\ref{ac0}), where $f$ satisfies {\bf (A1)}.
 Then, there exists a unique $c^*\in\R$ and a unique (up to a translation) traveling profile $U$ with
 speed $c^*$.
\end{theorem}
For a detailed proof, see Theorem~4.15 in \cite{F}.  It can be shown that 
the traveling profile $U$ is monotone increasing.
A phase portrait of the system (\ref{T2}) for various values of $c$
is sketched in Fig.~\ref{f:co5}.

For any given speed $c>c^*$, we seek a control in feedback form 
$\beta= \beta(u)\geq 0$, with finite cost, that yields a traveling wave with 
speed $c$.   % In the next section, we shall prove the existence of a control with minimum cost.
The main result of this section is

\begin{figure}[ht]
\centerline{\hbox{\includegraphics[width=8cm]{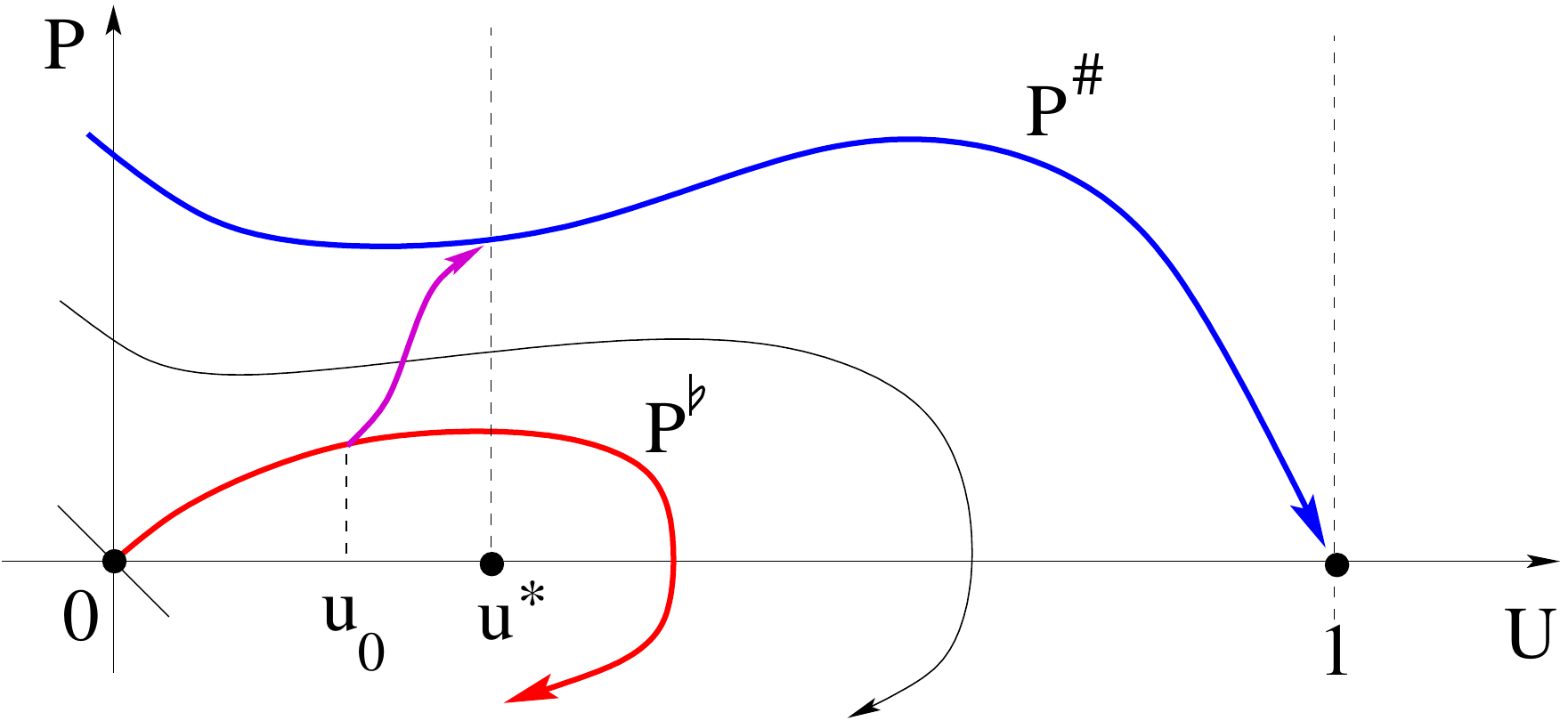}}}
\caption{\small Trajectories  of \eqref{ode3}  in the case $c>c^*$, $\beta(U)\equiv 0$. Here $P^\flat$ and $P^\sharp$
are the trajectories through $(0,0)$ and through $(1,0)$, respectively.
}
\label{f:co37}
\end{figure}

\begin{theorem}\label{t:22} Let $f$ satisfy the assumptions {\bf (A1)}  and let $c^*$ be as in
Theorem~\ref{t:21}. Then, for every $c>c^*$, there exist a bounded function $\beta:\,]0,1[\,\mapsto\R_+$
with compact support, such that the equation
\bel{tvc}
U'' + c U' + f(U) -\beta(U)~=~0,\qquad\qquad  U(-\infty) =0,\quad U(+\infty) =1.\eeq
admits a solution.
\end{theorem}
{\bf Proof.} {\bf 1.}
We will construct a solution of the first order system
\bel{ode3}
\left\{ \bega{rl} U'&= ~P,\\[1mm] 
P'&=~- c P - f(U) + \beta(U),\enda\right.\eeq
with asymptotic conditions
\bel{ac1}(U,P)(-\infty) \,=\,(0,0), \qquad\qquad (U,P)(+\infty)\,=\, (1,0),\eeq
for some function $\beta(\cdot)$ 
of the form
\bel{bdef}\beta(U)~=~\left\{ \bega{cl}  \gamma&\qquad \hbox{if}~~ u_0< U<u^*,\\[3mm] 0 &\qquad \hbox{otherwise.}\enda\right.
\eeq
 Here $u^*$ is the zero of $f$ considered 
in {\bf (A1)}, while $u_0\in \,]0,u^*[\,$ and $\gamma>0$ are suitable constants.
\v
{\bf 2.}  If $\beta\equiv 0$,
% by a linearization (see Chapter 4.4 in
 %\cite{F} for details) one checks that (0,0) is a node, while (1,0) is a saddle point.
computing the Jacobian matrix at a point $(U,P)$ one finds
\bel{JM}A(U,P)~=~\begin{pmatrix}
0 && 1\\[1mm]
-f'(U) && -c\end{pmatrix}.\eeq
Solving
$$\lambda^2 + c\lambda + f'(U)~=~0,$$
one obtains
\bel{la12}\lambda~=~{-c\pm\sqrt{c^2-4f'(U)}\over 2}\,.\eeq
We observe that the assumptions (\ref{f2})   imply that both $(0,0)$  and $(1,0)$
are saddle points.   In particular, the ODE
\bel{ODEc}
{d\over dU} P(U)~=~- c - {f(U)\over P}\eeq
has a solution $U\mapsto P^\flat(U)$ through $(0,0)$ with slope
$${dP^\flat \over dU}(0)~=~{-c+\sqrt{c^2-4f'(0)}\over 2}~>~0\,.$$
It also has a second solution $P^\sharp$ through the point $(1,0)$, with slope
$${dP^\sharp\over dU}(1)~=~{-c-\sqrt{c^2-4f'(1)}\over 2}~<~0\,.$$
In the special case where $c=c^*$, these solutions exactly match, as in Fig.~\ref{f:co5}, center.
On the other hand, when $c>c^*$, as shown in Fig.~\ref{f:co37} these two solutions satisfy
$$P^\flat(U)~<~P^\sharp(U) \qquad\qquad \forall U\in [0,u^*].$$
Now consider the backward Cauchy problem
\bel{ODE6}
{d\over dU} P(U)~=~- c - {f(U)\over P}+ \gamma,\qquad\qquad U\in [0, u^*],\eeq
with terminal data
\bel{td6}P(u^*)~=~P^\sharp(u^*).\eeq
By choosing $\gamma>0$ suitably large, the solution to (\ref{ODE6})-(\ref{td6}) will satisfy
$$P(U)< P^\flat(U)$$
at some point $0<U<u^*$.  Calling $u_0\in [0, u^*]$ the point where $P(u_0)=P^\flat (u_0)$, and defining $\beta(\cdot)$ as in
(\ref{bdef}), we  achieve the desired conclusion.
\endproof

\subsection{Existence of a control with finite cost.}

According to Theorem~\ref{t:22}, for every speed $c\geq c^*$ one can find a control $\beta=\beta(U)$
which yields a traveling wave with speed $c$.   However, in some cases such as (\ref{Lex}), one has
\bel{binf}\left\{\bega{rl}L(U,\beta)\,<\,+\infty \qquad \hbox{if}~\beta<\Hat \beta(U),\\[2mm]
L(U,\beta)\,=\,+\infty \qquad \hbox{if}~\beta\geq\Hat \beta(U),\enda\right.\eeq
for some function $\Hat \beta$. Therefore, some of the traveling waves considered in the above theorem may have infinite cost.

To understand in which cases a traveling wave exists  with finite cost,
consider any function $\Hat f$ that satisfies the assumptions on  $f$ stated in {\bf (A1)},
together with
\bel{fsh} f(u) - \Hat \beta(u)~\leq~\Hat f(u)~\leq~f(u)\qquad\qquad\forall u\in [0,1].\eeq
Call $\Hat c$ the speed of a traveling wave for the corresponding equation
$$u_t~=~u_{xx} + \Hat f(u).$$

\begin{theorem}\label{t:23}  In the above setting, for every speed $c\in [c^*, \Hat c[\,$ 
there exists a control $\beta=\beta(u)$
with finite cost,  such that the equation (\ref{tvc}) has a solution.
%Namely, there exists a traveling wave with speed $c$.
\end{theorem}

\begin{figure}[ht]
\centerline{\hbox{\includegraphics[width=11cm]{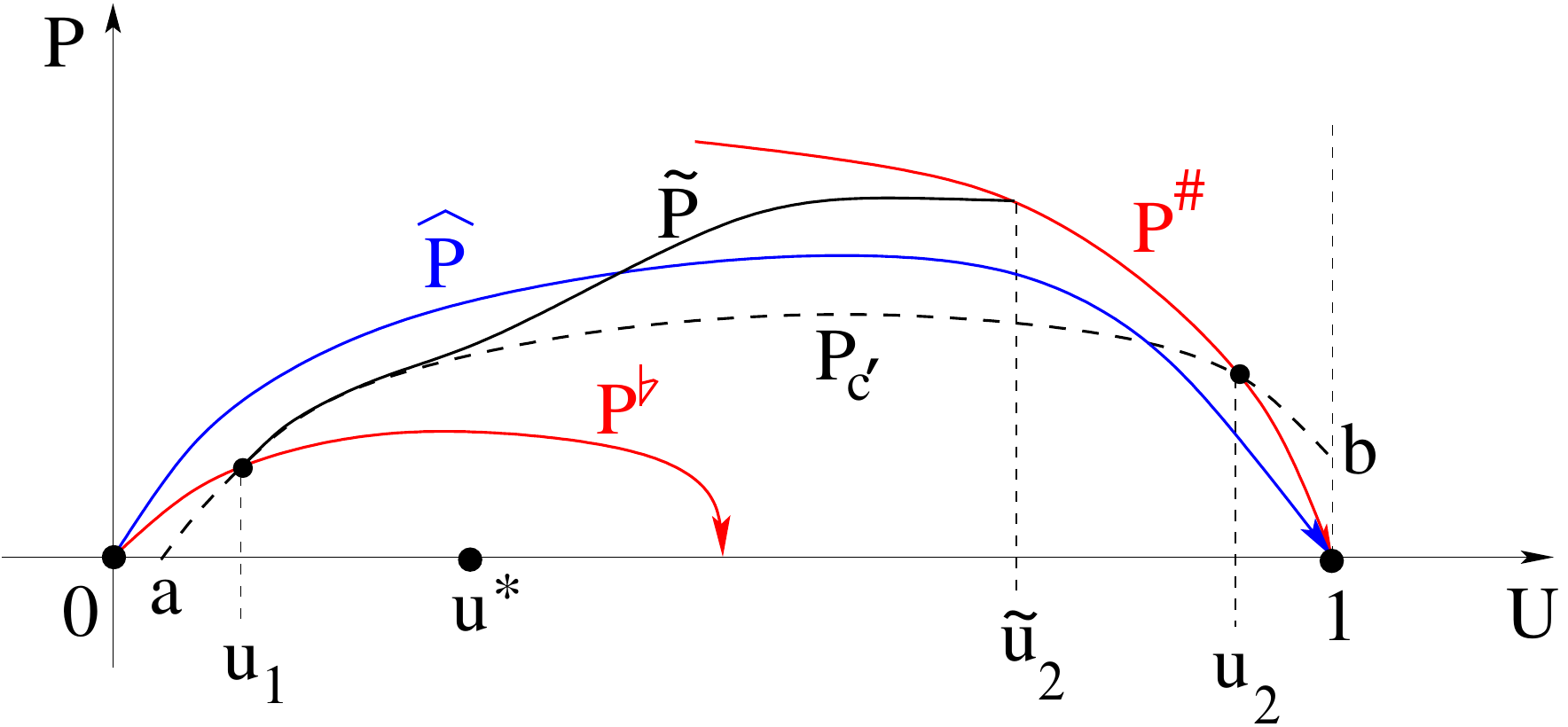}}}
\caption{\small The trajectories considered in the proof of Theorem~\ref{t:23}.}
\label{f:co43}
\end{figure}

{\bf Proof.}   {\bf 1.} We can assume $\Hat c>c^*$, since otherwise there is nothing to prove.
By assumption,
the system
\bel{tp3}
\left\{ \bega{rl} U'&= ~P,\\[1mm] 
P'&=~- \Hat c P - \Hat f (U) \enda\right.\eeq
has a heteroclinic orbit joining (0,0) with (1,0).
With reference to Fig.~\ref{f:co43}, we call $P=\Hat P(U)$ the corresponding solution to
$${dP\over dU}~=~- \Hat c -{ \Hat f(U)\over P}\,.$$
In addition, we denote by $P=P^\flat(U)$  and $P=P^\sharp(U)$ the solutions to 
$${dP\over dU}~=~- c -{ f(U)\over P}\,,$$
with boundary data
$$P^\flat (0)\,=\,0\qquad\hbox{and}\qquad P^\sharp(1)\,=\,0,$$
respectively.
\v
{\bf 2.} Next, choose any speed $c'$ with $$c^*\,< \,c\,<\,c'\,<\,\Hat c.$$   
The strict inequality $c'<\Hat c$ implies that the system
\bel{tp30}
\left\{ \bega{rl} U'&= ~P,\\[1mm] 
P'&=~- c' P - \Hat f (U) ,\enda\right.\eeq
has an orbit joining a point $(a,0)$ on the positive $U$-axis  with a point
$(1,b)$, with $a,b>0$.   We call $P=P_{c'}(U)$ this profile.
\v
{\bf 3.} Still referring to Fig.~\ref{f:co43}, consider the intersection points
$0<u_1<u_2<1$, defined by
$$P^\flat(u_1)~=~P_{c'}(u_1),\qquad\qquad P^\sharp(u_2)~=~P_{c'}(u_2).$$
Define the control
\bel{tilu}\Tilde\beta(U)~\doteq~\max\bigl\{  \Hat \beta(U) -(c'-c) P_{c'}(U)\,,~0   \bigr\} .  \eeq
Notice that this implies
$$-c' - {f(U)-\Hat \beta(U)\over P_{c'}(U)}~\leq~-c - {f(U) -\Tilde\beta(U)\over P_{c'}(U)}
\qquad \forall U\in [u_1, u_2].$$
Calling $P=\Tilde P(U)$ the solution to
$${dP\over dU}~=~ - c -  {f(U) -\Tilde\beta(U)\over P}\,,\qquad\qquad \Tilde P(u_1) \,=\,
P_{c'}(u_1),$$
a comparison argument yields
\bel{TTU}\Tilde P(U)\geq P_{c'}(U)\qquad \forall U>u_1\,.\eeq
Therefore, the curve $P=\Tilde P(U)$ will intersect the trajectory $P^\sharp$ at some point
$\Tilde u_2\leq u_2$.
\v
{\bf 4.}  
We claim that the concatenation of trajectories
\bel{pco}
P(U)~=~\left\{\bega{rl} P^\flat(U)\quad&\hbox{if}\quad U\in [0, u_1],\\[1mm]
\Tilde P(U)\quad&\hbox{if}\quad U\in [u_1, \Tilde u_2],\\[1mm]
P^\sharp(U)\quad&\hbox{if}\quad U\in [ \Tilde u_2, 1],\enda\right.\eeq
provides a solution to (\ref{T2}) with finite cost.

Indeed,  for $U\in [0, u_1]\cup [\Tilde u_2,1]$ the above solution corresponds to a control
$\beta=0$, with zero cost.  

Furthermore, for $U\in [u_1, \Tilde u_2]$, in view of (\ref{cp1}) the cost is 
\bel{cost5}\int_{u_1}^{\Tilde u_2} {L\bigl( U, \Tilde \beta(U)\bigr)\over \Tilde P(U)}\, dU.\eeq
By (\ref{tilu}) we have
$$\Hat \beta(U) - \Tilde\beta(U)~>~\delta~>~0\qquad\qquad \forall U\in [u_1, \Tilde u_2].$$
Hence the numerator 
$L\bigl( U, \Tilde \beta(U)\bigr)$ remains uniformly bounded for $U\in [u_1, \Tilde u_2]$.
Finally, the denominator 
$\Tilde P(U)$ is uniformly positive, because of (\ref{TTU}).
\endproof

\section{Existence of an optimal strategy}
\label{sec:3}
\setcounter{equation}{0}

Extending one of the results in \cite{BCS1} to this more general nonlinear setting, we now prove

\begin{theorem}\label{t:31} Let $f,L$ satisfy the assumptions {\bf (A1)} and {\bf (A2)}.   Then,
for any wave speed $c>c^*$, if (\ref{cp2}) has a solution with finite cost $J(U,\beta)<\infty$, then the  problem {\bf (OTW)} 
has an optimal solution.
\end{theorem} 

{\bf Proof.} {\bf 1.} 
Following the direct method in the Calculus of Variations, we consider a minimizing sequence $(u_n, \beta_n)_{n\geq 1}\,$.
That is, a sequence of solutions to (\ref{cp2}) such that
\bel{min} \lim_{n\to \infty} \int_{-\infty}^{\infty} L(u_n, \beta_n)\, dx~=~\inf_{(u,\beta)} \int_{-\infty}^{\infty} L(u, \beta)\, dx.\eeq
Here the infimum is taken over all solutions $(u,\beta)$ of  (\ref{cp2}).
 By a translation in the $x$-variable, we can assume that 
\bel{un0} u_n(0)~=~u^*\qquad\forall n\geq 1.\eeq
\v
{\bf 2.} By  the growth condition (\ref{LL}), it follows that the  norms $\|\beta_n\|_{\L^p}$
are uniformly bounded.
\v
{\bf 3.} In this step we prove that the functions $u_n$ are uniformly Lipschitz continuous.
Since all these functions are nondecreasing, it suffices to show that their derivative $p_n(x)= u'(x)$ is 
bounded above, uniformly for all $x\in \R$.
%
%{\color{red} Hint: under the growth condition (\ref{LL}), we have
%$P\in W^{1,p}(\R)$, $U\in W^{2,p}(\R)$.   Hence we have a uniform bound on $P = U'$.}
%
Calling $M$ the maximum value of  the function $f$ on $[0,1]$, from (\ref{ode3}) it follows
\bel{p_n}p_n'(x)~=~-c\cdot p_n(x) -f(u_n(x))+\beta_n~\geq~-c p_n(x)- M. \eeq
In turn, for any $x_0\in \R$ this yields the lower bound
%Consider that $u_n$ and $p_n(x_0)$ for one $x_0$ is known and solve $p_n$ out on $[x_0,+\infty)$ as:
$$ p_n(x) = e^{-c(x-x_0)}p_n(x_0)-\int_{x_0}^x e^{-c(x-\xi)}M \,d\xi.$$
Integrating the above equation above from $x_0$ to $x_0+1$, and observing that $u_n(x)\in [0,1]$, in the case $c\not= 0$
we obtain
$$\bega{rl} 1&\ds \geq~u_n(x_0+1)-u_n(x_0) ~=~\int_{x_0}^{x_0+1} p_n(x)\, dx\\[4mm]
&\ds \geq~\int_{x_0}^{x_0+1} e^{-c(x-x_0)}p_n(x_0)\, dx-
\int_{x_0}^{x_0+1} \int_{x_0}^x M \,e^{-c(x-\xi)}d\xi\,dx\\[4mm]
&\ds =~  {1-e^{-c}\over c}p_n(x_0) -M \cdot \Big({1\over c}+{ e^{-c}-1\over c^2}\Big).
 \enda$$
This yields the bound
 $$p_n(x_0)~\leq~ {c \over 1-e^{-c}}+M\Big({1 \over 1-e^{-c}}-{1\over c}\Big).$$
 Notice that this bound is uniformly valid for every $x_0\in \R$ and $n\geq 1$. 
We thus conclude that all functions $u_n$ have uniformly bounded derivatives, hence are uniformly Lipschitz continuous.

In the case $c=0$, the above computation is simply replaced by
$1\geq p_n(x_0) - {M\over 2}$, leading to the same conclusion.

\v
{\bf 4.} Since all functions $u_n$ are uniformly Lipschitz continuous, by possibly
taking a subsequence, we can assume the convergence
\bel{unc}u_n(x)~\to~u(x)\eeq
uniformly for $x$ in bounded sets.
Moreover, since the $\L^p$ norms of the functions $\beta_n$ are uniformly bounded, 
we have the weak convergence $\beta_n\wto \beta$ for some function
$\beta\in \L^p(\R)$.

We can write the differential equation satisfied by $u_n$ in integral form:
\bel{uni}
u_n'(x_2)-u_n'(x_1) + cu_n(x_2)-cu_n(x_1) = \int_{x_1}^{x_2} \big[-f(u_n) + \beta_n \big] \,dx,
\eeq
which is valid for every $x_1<x_2$.  Taking the limit as $n\to\infty$ and recalling the uniform convergence $u_n\to u$
and the weak convergence $\beta_n\wto \beta$ in $\L^p$, we obtain
\bel{upl}
u'(x_2)-u'(x_1) + cu(x_2)-cu(x_1) = \int_{x_1}^{x_2} \big[-f(u) + \beta \big] \, dx.\eeq
Therefore,  $u\in W^{2,p}_{loc}(\R)$, and  the ODE in \eqref{cp2} is satisfied.

\v
{\bf 5.} In the next two steps, using the assumptions {\bf (A2)} on $L$, we prove the lower semicontinuity relation:
\bel{LSC} 
 \int_{-\infty}^{+\infty} L\bigl(u(x), \beta(x)\bigr)\, dx~\leq~
 \liminf_{n\to \infty}~ \int_{-\infty}^{+\infty} L\bigl(u_n(x), \beta_n(x)\bigr)\, dx.\eeq
Since $L\geq 0$, we have
\bel{LSC1}  \int_{-\infty}^{+\infty} L\bigl(u(x), \beta(x)\bigr)\, dx~=~\lim_{R\to +\infty} 
 \int_{-R}^R L\bigl(u(x), \beta(x)\bigr)\, dx\,.\eeq
 Next, for every $m\geq 1$ consider the function
 \bel{Lm}L^{(m)}(u,\beta)~\doteq~\min_{z\in [0,\beta]} \Big\{ L(u,z) +m(\beta-z)\Big\}.\eeq
 In view of {\bf (A2)}, $L^{(m)}$ is continuous w.r.t.~both variables $u,\beta$, and  Lipschitz continuous
 with constant $m$ in the variable $\beta$, uniformly for every $u$.  Indeed, $L^{(m)}$ is the largest function $\leq L$ with these properties.
 Since
 $$L(u,\beta)~=~\lim_{m\to\infty}  L^{(m)}(u,\beta),$$
 we have
 \bel{Lm3}\int_{-R}^R L(u(x), \beta(x)\bigr)\, dx~=~\lim_{m\to\infty}
 \int_{-R}^R L^{(m)}\bigl(u(x), \beta(x)\bigr)\, dx.\eeq
To prove (\ref{LSC}), it thus suffices to show that 
\bel{LSC4} 
 \int_{-R}^R L^{(m)}\bigl(u(x), \beta(x)\bigr)\, dx~\leq~
 \liminf_{n\to \infty}~ \int_{-R}^R L^{(m)}\bigl(u_n(x), \beta_n(x)\bigr)\, dx,\eeq
for any given $R, m\geq 1$.
\v
{\bf 6.} By the
convexity of the maps $\beta\mapsto L^{(m)} \bigl(u(x),\beta\bigr)$ it follows
%for every $x\in [-R,R]$ we have
%\bel{Lun}L^{(m)}(u(x),\beta_n(x))~\geq~L^{(m)}(u_n(x), \beta(x)) + L^{(m)}_\beta(u_n(x), \beta(x))\cdot (\beta_n(x)- \beta(x)). \eeq
%{\color{red} NOTE: we already know that $L^{(m)}$ is uniformly Lipschitz continuous in $\beta$.
%No need to introduce the set $A_M$.}
\bel{iam}\bega{l} \ds\int_{-R}^R L^{(m)}\bigl(u_n(x),\beta_n(x)\bigr)\, dx\\[4mm]
\qquad\ds \geq~ 
\int_{-R}^R L^{(m)}\bigl(u_n(x), \beta(x)\bigr) \,dx + \int_{-R}^R  L^{(m)}_\beta\bigl(u_n(x), \beta(x)\bigr)\cdot \bigl(\beta_n(x)- \beta(x)\bigr)\, dx.\enda\eeq
Thanks to the uniform bound
$$\bigl| L^{(m)}(u,\beta)\bigr|~\leq~m \beta\qquad\qquad\forall u\in [0,1], ~~\beta\geq 0,$$
and the uniform convergence  $u_n\to u$, the first integral on the right hand side of 
(\ref{iam}) satisfies
\bel{liam1} \lim_{n\rightarrow \infty} \int_{-R}^R L^{(m)}\bigl(u_n(x), \beta(x)\bigr)\, dx~=~
 \int_{-R}^R L^{(m)}\bigl(u(x), \beta(x)\bigr)\, dx\,.\eeq
Using the uniform convergence  $u_n\to u$, the weak convergence $\beta_n\wto\beta$,
and
observing that $\|L_\beta\|_{\L^\infty}\leq m$,
we conclude that the second integral on the right hand side of 
(\ref{iam}) satisfies
\bel{liam2} \lim_{n\rightarrow \infty} \int_{-R}^R L^{(m)}_{\beta}\bigl(u_n(x), \beta(x)\bigr)\cdot\bigl(\beta_n(x)-\beta(x)\bigr)\,dx \,= \,0. \eeq
Together, (\ref{liam1})-(\ref{liam2}) imply (\ref{LSC4}), and hence (\ref{LSC}).
\v
\begin{figure}[ht]
\centerline{\hbox{\includegraphics[width=9cm]{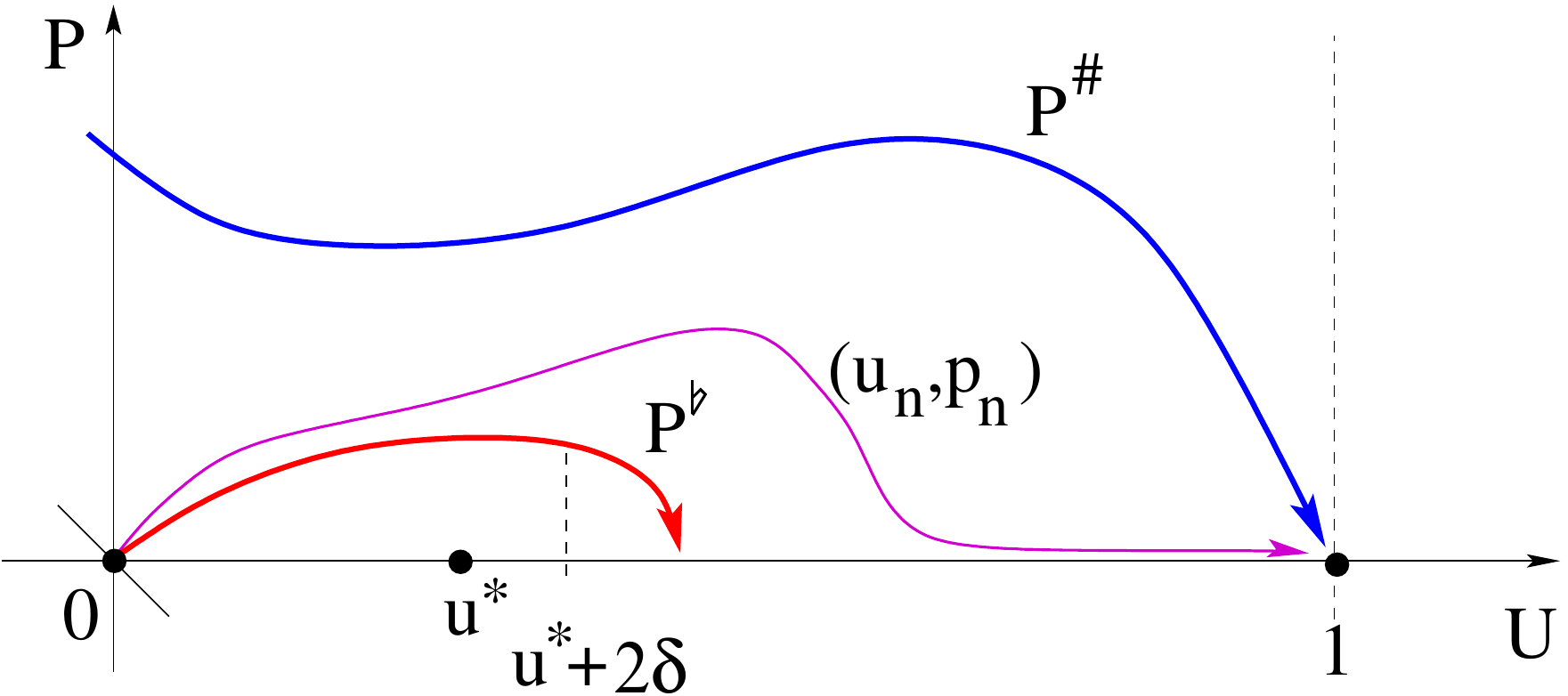}}}
\caption{\small All trajectories $x\mapsto \bigl(u_n(x), p_n(x)\bigr)$ take values in the region
between $P^\flat$ and $P^\sharp$.
}
\label{f:co38}
\end{figure}

{\bf 7.}  In this step we complete the proof by establishing the limits
\bel{limu}
u^-\,\doteq\,\lim_{x\to -\infty} u(x)\,=\,0,\qquad\qquad  u^+\,\doteq\,\lim_{x\to +\infty} u(x)\,=\,1.\eeq
Notice that the fact that every $u_n$ satisfies the above limits, together with the convergence
$u_n\to u$ uniformly on bounded sets, does not suffice to conclude (\ref{limu}).   
From the monotonicity of $u$ it only follows that the limits $u^-, u^+$ exists, with $0\leq u^-\leq u^+\leq 1$.
To achieve (\ref{limu}), a more careful argument is needed.

With reference to Fig.~\ref{f:co38}, we observe that any trajectory $x\mapsto \bigl(u_n(x), p_n(x)\bigr)$ is contained in the region between $P^\flat$ and $P^\sharp$, i.e. the unstable manifold  through $(0,0)$ and the stable manifold through $(1,0)$, respectively.    In particular, for some constants  $C,\delta>0$ independent of $n$, we have the implication
\bel{pnla} u_n(x)\in [0, u^*+2\delta]\quad\implies\qquad u'_n(x)~=~p_n(x)~\geq ~C u_n(x).\eeq
This immediately implies 
\bel{umz}u^-\,=\,0,\qquad\qquad u^+\,\geq\,u^*+2\delta.\eeq
It remains to prove that $u^+=1$.
 We prove this claim by contradiction. 
By (\ref{umz}) we can find  $R>0$ large enough so that  $u(x)>u^*+\delta$ for all  $x>R$. 
Integrating the differential equation in \eqref{cp2} on the interval  $[R, 2R]$, we obtain
\bel{intb}
u'(2R)-u'(R) + c\bigl(u(2R)-u(R)\bigr) + \int_R^{2R} f(u(x))dx~ =~ \int_R^{2R} \beta(u(x))dx.\eeq
Since $u$ is Lipschitz continuous, the quantities 
$$u'(2R),\quad u'(R),\quad u(2R),\quad u(R),$$ are all uniformly bounded.  On the other hand,
\bel{minf}\int_R^{2R} f(u(x))dx~\geq ~R\cdot \min_{u\in [u^*+\delta,\,  u^+]} f(u).\eeq
 If $u^+<1$, then the minimum in (\ref{minf}) is strictly positive.  Hence the right hand side of (\ref{minf}) approaches infinity as 
 $R\to +\infty$.  In particular, choosing $R$ large enough, from (\ref{intb}) we obtain
$$ \int_R^{2R} \beta(u(x))dx~\geq~1.$$
Repeating the same argument on the intervals $\bigl[k R, (k+1) R\bigr]$, we again obtain
\bel{ikr} \int_{kR}^{(k+1)R} \beta(u(x))dx~\geq~1,\qquad\qquad k=1,2,3,\ldots\eeq
By the assumption (\ref{LL}), in view of (\ref{ikr}), this implies
\bel{IU}\bega{l} J(u,\beta)~\ds\geq~\int_R^{+\infty} L\bigl(u(x),\beta(x)\bigr)\, dx~=~\sum_{k\geq 1} 
\int_{kR}^{(k+1)R} L\bigl(u(x),\beta(x)\bigr)\, dx\\[4mm]
\qquad \ds \geq ~\sum_{k\geq 1} 
\int_{kR}^{(k+1)R} C_1\beta^p(x)\, dx~ \geq ~\sum_{k\geq 1} 
\int_{kR}^{(k+1)R} C_1\left({1\over R}\right)^p dx ~=~\sum_{k\geq 1}  C_1 R^{1-p}~=~+\infty. 
\enda\eeq
We thus obtain a contradiction with the previous step, where we proved that the cost $J(u,\beta)$ is finite.
This completes the proof.
\endproof

\section{Necessary conditions for optimality}
\label{sec:4}
\setcounter{equation}{0}
Given a speed $c>c^*$,
assume that $(U,\beta)$ yield an optimally controlled traveling wave profile, as in Theorem~\ref{t:31}.
We seek necessary conditions to determine this profile.

In terms of the  $U$-$P$ coordinates, as in (\ref{ode3}), this means that the control 
$\beta= \beta(U)\geq 0$ minimizes the cost functional
\bel{cost1}  J(\beta)~=~
\int_0^1 {L\bigl( U,  \beta(U)\bigr)\over P (U)}\, dU,\eeq
subject to
\bel{PUU}
{dP\over dU}~=~- c + {\beta-f(U)\over P} ,\qquad P(0)= P(1)=0.\eeq
To apply the Pontryagin Maximum Principle \cite{BPi, Ce},
%Here $U\mapsto \beta(U)\geq 0$ is any measurable control function, subject only to the 
%non-negativity condition. 
%This is a typical problem of optimal control with initial and terminal constraints, in Bolza form.
we first compute
$${\partial\over\partial P} \left({L( U,  \beta)\over P}\right)~=~- {L( U,  \beta)\over P^2}\,,\qquad\qquad 
{\partial\over\partial P} \left(- c + {\beta-f(U)\over P} \right)~=~ -{\beta-f(U)\over P^2}\,.$$
The PMP now yields the existence of an adjoint variable $Y(\cdot)$ satisfying
the linear equation
\bel{Yeq} 
{dY\over dU} ~=~ {\beta(U)-f(U)\over P^2(U)}\,Y+ {L\bigl( U,  \beta(U)\bigr)\over P^2(U)} \,,
\eeq
such that, at a.e.~$U\in [0,1]$, the following optimality condition holds: 
\bel{min11} \beta(U)~=~\argmin_{\beta\geq 0} \left\{ \left( -c+{\beta-f(U)\over P(U)} \right)Y(U)  +  {L\bigl( U,  \beta\bigr)\over P (U)}\right\}.
\eeq
Equivalently,
\bel{min22}\beta(U)~=~\argmin_{\beta\geq 0} \Big\{ \beta Y(U) + L (U,\beta)\Big\}.
\eeq

Note that, in the region where $\beta(U)>0$,  by (\ref{min22}) we must have
\bel{min3} Y(U)+ L_\beta\bigl( U, \beta(U)\bigr)~=~0.\eeq
Differentiating w.r.t.~$U$ and using (\ref{Yeq})-(\ref{min3}), we obtain
\bel{YU} \bega{l}\ds {d\over dU}  L_\beta\bigl( U, \beta(U)\bigr) +{d\over dU} Y(U)~=~
\ds{d\over dU}  L_\beta\bigl( U, \beta(U)\bigr) +
 {\beta(U)-f(U)\over P^2(U)}\,Y+ {L\bigl( U,  \beta(U)\bigr)
\over P^2(U)}\\[4mm]
\qquad \ds=~\ds{d\over dU}  L_\beta\bigl( U, \beta(U)\bigr) -
 {\beta(U)-f(U)\over P^2(U)}\, L_\beta\bigl( U, \beta(U)\bigr)+ {L\bigl( U,  \beta(U)\bigr)
\over P^2(U)}~=~
0.
\enda
\eeq
In most cases, the control  $\beta$ will be active only on some interval $]u_1, u_2[\,$,  so that
\bel{beta}
\left\{\bega{rl}\beta(U) >0\quad &\hbox{if}~~U\in\,]u_1,u_2[\,,\\[2mm]
\beta(U)=0\quad &\hbox{if} ~~U\in  [0,u_1]\cup [u_2, 1]\,.\enda\right.\eeq
From (\ref{min22}) and the strict convexity of $L(u,\cdot)$ it now follows
$$\lim_{U\to u_1+} \beta(U)~=~\lim_{U\to u_2-} \beta(U)~=~0.$$

The optimal solution can thus be obtained by solving the  ODE  in (\ref{YU})
over an interval $[u_1, u_2]$, whose endpoints are determined by 
the two additional boundary 
conditions
\bel{BC12} 
Y(u_1) + L_\beta(u_1,0)~=~Y(u_2) + L_\beta(u_2,0)~=~0.\eeq

\begin{figure}[ht]
\centerline{\hbox{\includegraphics[width=16cm]{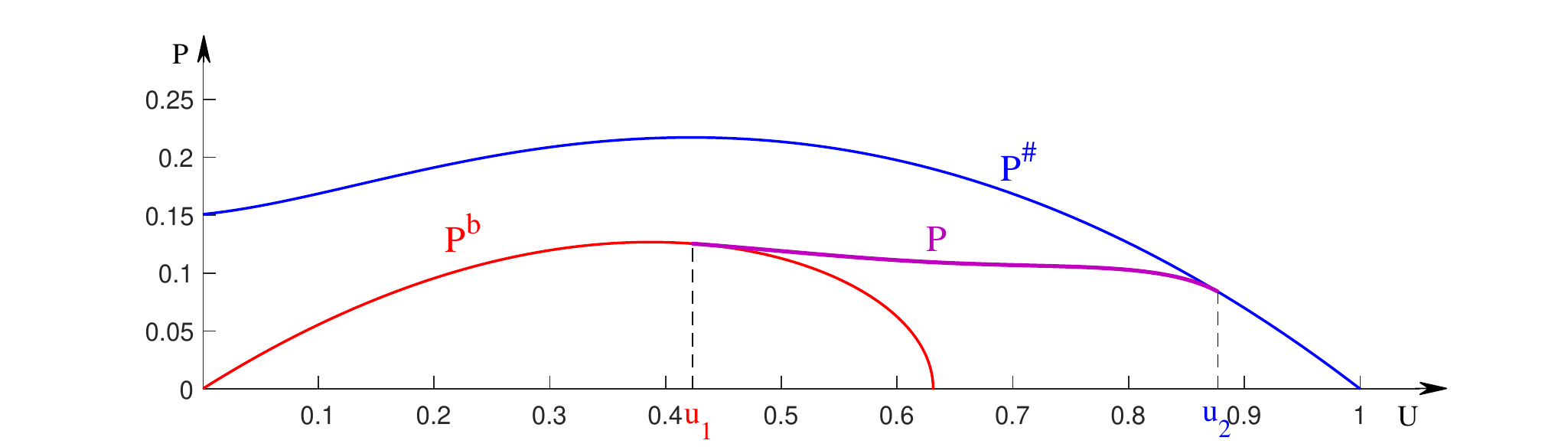}}}
\caption{\small   The optimal traveling profile for the given speed $c=-0.1$, in the $U,P$ coordinates.}
\label{f:c01profile}
\end{figure}

\begin{figure}[ht]
\centerline{\hbox{\includegraphics[width=16cm]{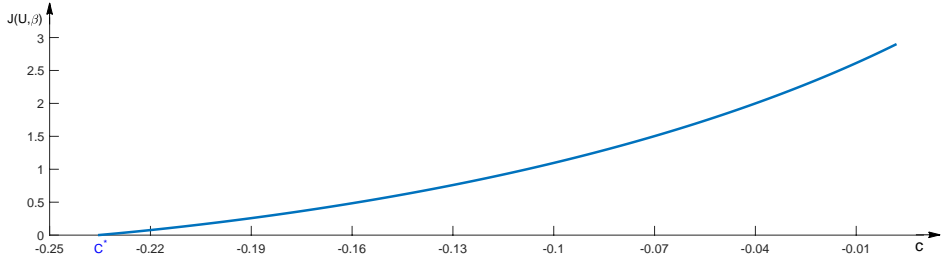}}}
\caption{\small The minimum cost $E(c)$, depending on the wave speed $c\geq c^*$.}
\label{f:minco}
\end{figure}

 \subsection{Numerical computation of optimally controlled traveling profiles.}

In a typical application, the optimal traveling profile with a given speed $c>c^*$ can be 
computed as follows.
\v
STEP 1: Observing that both $(0,0)$ and $(1,0)$ are both saddle points for the system (\ref{T2}), compute the unstable manifold $P=P^\flat(U)$ through $(0,0)$, and the 
stable manifold $P=P^\sharp(U)$ through $(1,0)$,
as shown in Fig.~\ref{f:c01profile}.
\v
STEP 2: Determine the interval $[u_1, u_2]$ and the portion of the optimal trajectory $U\mapsto P(U)$ 
for $U\in [u_1, u_2]$ by solving the system of two equations
\bel{ELS}\left\{
\bega{rl} \ds
{dP\over dU}&\ds=~- c + {\beta-f(U)\over P} \,,\\[4mm]
\ds{d\beta\over dU} & \ds=~{1\over L_{\beta\beta} ( U, \beta)}
\cdot\left[  {\beta-f(U)\over P^2(U)}\, L_\beta(U,\beta) - {L( U,  \beta)
\over P^2(U)} - L_{U\beta }(U,\beta)\right],\enda\right.\eeq
with the four boundary conditions
\bel{BC4}\left\{\bega{rl} P(u_1) &=~P^\flat (u_1),\\[2mm]
 P(u_2) &=~P^\sharp (u_2),\enda\right.\qquad\qquad \left\{\bega{rl}\beta(u_1) &=~0,\\[2mm]
 \beta(u_2) &=~0.\enda\right.
 \eeq
Note that the solution to a system of two first order ODEs is determined by two
boundary conditions.  Here the two additional conditions $\beta(u_1)=\beta(u_2)=0$ are needed 
to determine the endpoints of the interval $[u_1,u_2]$.
\v
\begin{example} {\rm
For sake of illustration, 
we consider here the optimization problem for a traveling profile, 
choosing  $f(u)$ and $L(u,\beta)$ as in (\ref{fu3})-(\ref{Lex}), with $u^* = 1/3$. 
If no control is present, a numerical simulation shows that the speed of the traveling profile
is 
$c^*\approx -0.2356$.

For various speeds $c>c^*$, we seek the minimum cost of a control that produces a traveling profile with speed $c$.
This is achieved following the above steps 1 and 2.
To achieve Step 2,
for a given $u_1$, consider the solution $U\mapsto \bigl(P(U), \beta(U)\bigr)$ of 
(\ref{ELS}) with initial data
\bel{BC6} P(u_1) ~=~P^\flat(u_1), \qquad\qquad \beta(u_1)~=~0.\eeq
This solution is prolonged until $P(U) = P^\sharp (U)$.  More precisely, 
let $u_2>u_1$ be the first point such that
$$P(u_2)~=~P^\sharp (u_2).$$
The above construction yields a map 
$$u_1~\mapsto~ \phi(u_1) ~\doteq~\beta(u_2).$$
By a shooting method, we determine $u_1\in [u^*, 1]$ such that $\phi(u_1)=0$. This yields the desired solution.

In the case where the wave speed is $c=-0.1$, a numerical simulation of the optimal traveling profile,
is shown in Fig.~\ref{f:c01profile}.
The minimum cost, for increasing values of the wave speed $c\geq c^*$, is shown in Fig.~\ref{f:minco}.
}\end{example}

\section{Controlled traveling profiles for Model 1}
\label{sec:5}
\setcounter{equation}{0}
In this section we consider controlled traveling profiles for the system (\ref{ode1}), say with
$$u(t,x)\,=\, U(x-ct),\qquad\quad \theta(t,x)\,=\, \Theta(x-ct),\qquad\quad \alpha=\alpha(x-ct).$$
Since (\ref{ode1}) is in triangular form,   for any speed $c>c^*$ the existence
of an optimal  traveling profile $U$ for the first equation has already 
been proved in Theorem~\ref{t:31}.
The next result shows that, if $c<0$, then the second equation in (\ref{ode1}) also admits
a traveling profile with speed $c$. 

% \bel{ode1} \left\{\bega{rl} u_t&=~u_{xx} + f(u) - \beta,
% \\[1mm]
% \theta_t&=~\kappa_2\,u\,(1-\theta).\enda\right.
% \eeq
We recall that the functions $x\mapsto \bigl( U(x), \Theta(x)\bigr)$  should  satisfy
\bel{tv3}U'' + c U' + f(U, \alpha(U))~=~0, \qquad\qquad U(-\infty) =0,\quad U(+\infty)=1,\eeq
\bel{tv4} c\Theta' + \kappa_1 U (1-\Theta)~=~0,\qquad \qquad \Theta(-\infty) = 0, \qquad\Theta(+\infty) =1.\eeq
A solution $\Theta$  of (\ref{tv4}) will be constructed assuming the integrability condition
\bel{Uint}
\int_{-\infty}^0 U(x)\, dx~ <~+\infty.\eeq

\begin{theorem}\label{t:51}
Let $U:\R\mapsto [0,1]$ be an increasing  solution to (\ref{tv3}), such that  (\ref{Uint}) holds.
 Then a  solution to (\ref{tv4}) exists  if and only if
$c<0$.\end{theorem}

{\bf Proof.}  To construct the function $\Theta$ in (\ref{tv4}), we begin by solving
$${-\Theta' \over 1-\Theta}~=~{\kappa_1\over c}U, \qquad \qquad \Theta(-\infty)=0.$$
%We suppose that $$ always holds true. In this way, only the limit of $\Theta$ at $+\infty$ should be paid attention to. So we integrate the equation above from $-\infty$ to $x<+\infty$,
An integration yields
$$\ln\big(1-\Theta(x)\big)~=~{\kappa_1\over c}\int_{-\infty}^x U(y)\,dy,$$
%Notice that the integrability assumption on $U$ implies
$$\Theta(x)~=~1 - \exp\left\{ {\kappa_1\over c}\int_{-\infty}^x U(y)\,dy\right\}.$$ 
Since $\kappa_1>0$, if  $c<0$ then
$$\lim_{x\to -\infty} ~\Theta(x)~=~0, \qquad\qquad \lim_{x\to +\infty} ~\Theta(x)~=~1.$$
On the other hand, if $c>0$ then 
$$\lim_{x\to +\infty} ~\ln\big(1-\Theta(x)\big)~=~\lim_{x\to +\infty}{\kappa_1\over c}\int_{-\infty}^x U(y)\,dy ~=~+\infty.$$
This contradicts the condition $\Theta(x)\in [0,1]$. 
Hence, no such traveling profile exists.\endproof

\v
A key assumption of the previous theorem was the boundedness of the integral in  (\ref{Uint}). 
We now show that this is always satisfied in the setting considered in Theorems~\ref{t:22} and \ref{t:31}.

\begin{lemma}\label{l:42} Assume that $f:[0,1]\mapsto\R$ satisfies the assumptions {\bf (A1)}.   
Then for any $c> c^*$ and any solution $U$ of (\ref{tvc}) with $\beta(U)\geq 0$, the integrability condition
(\ref{Uint}) holds.
\end{lemma} 

{\bf Proof.} As remarked in (\ref{pnla}), under the assumptions {\bf (A1)} any traveling wave solution must satisfy
\bel{U'} U'(x)~\geq~C \, U(x)\qquad\quad \hbox{whenever} \qquad U(x)\in [0, u^*+2\delta],\eeq
for some positive constants $C, \delta>0$. 
Calling $x^*\in \R$ the point where $U(x^*) = u^*$, from the  differential inequality (\ref{U'}) we deduce
\bel{uex}
U(x)~\leq~e^{-C(x^*-x)} u^* \qquad\qquad\forall x\in \,]-\infty, x^*].
\eeq
This implies that, as $x\to -\infty$, the function $U(x)$ converges to zero exponentially fast.
Hence the integrability condition (\ref{Uint}) holds.\endproof

\section{Traveling profiles for Model 2}
\label{sec:6}
\setcounter{equation}{0}
In this section we begin a study of the system (\ref{ode2}), assuming that
the function $f$ satisfies the assumptions in {\bf (A1)} together with
\bel{drate}
f(u)~\geq~-d u\qquad\qquad u\in [0,1].\eeq
Introducing the variable
$v= Iu =$ density of infected insects, we thus consider the system
\bel{71}
\left\{ \bega{rl} u_t&=~ u_{xx} + f(u)- \alpha u\,,\\[2mm]
v_t &=~ v_{xx} + \kappa_2 (u-v) \theta - \alpha v - d \,v ,\\[2mm]
\theta_t&=~ \kappa_1  (1-\theta) \,v\,.\enda\right.\eeq

For future use, we recall a basic definition~\cite{F,VVV}.

\begin{definition}   A $\C^1$ function $F:\R^m\mapsto\R^m$, 
say $F(w)= \bigl(F_1(w),\ldots, F_m(w)\bigr)$ is {\bf quasi-monotone} on a convex 
domain $\D\subseteq\R^m$ if
$${\partial F_i\over\partial w_j}(w)~\geq~0\qquad\qquad\forall i\not= j,~~~w=(w_1,\ldots, w_m)\in \D.$$
\end{definition}
Motivated by (\ref{71}), we observe that the map $F:\R^3\mapsto\R^3$ defined by
\bel{Fdef}F(u,v,\theta)~=~\Big( f(u)-\alpha(x) u\,,~ \kappa_2 (u-v) \theta - \alpha(x) v - d \,v \,,~
\kappa_1(1-\theta) v\Big),\eeq
 is quasi-monotone on the 
domain
\bel{PID}\D~\doteq~\Big\{ (u,v,\theta)\,;\qquad 0\leq v\leq u\leq 1,\quad \theta\in [0,1]\Big\}.\eeq

By a comparison argument we obtain
\begin{lemma}\label{l:61} Let $f$ satisfy the assumptions {\bf (A1)} together with the inequality (\ref{drate}).  
Then the domain $\D$ is positively invariant for  the system (\ref{71}). Namely, for any control
function $\alpha=\alpha(t,x)\geq 0$,
let $(u,v,\theta)$ be
a solution to (\ref{71}) such that, at time $t=0$,
 $(u,v,\theta)(0,x)\in\D$ for all $x\in\R$. 
 Then  $(u,v,\theta)(t,x)\in\D$ for all $x\in\R$ and $t\geq 0$. 
\end{lemma}

{\bf Proof.} We first observe that the triples 
$$(u^-, v^-, \theta^-)(t,x)\,=\, (0,0,0),\qquad\qquad (u^+, v^+, \theta^+)(t,x)\,=\, (1,1,1),$$ 
provide a subsolution and a supersolution to the system (\ref{71}), respectively.   
This implies that the three functions $u,v,\theta$ all take values within the interval $[0,1]$.

Next, let $(u,v,\theta)$ be any solution.  Then  the function $w= u-v$ satisfies  
\bel{wpso}
w_t~=~u_t - v_t~=~\Delta w + \bigl[ f(u) +dv\bigr] -\kappa_2 \theta w -\alpha w~\geq~\Delta w -
\bigl[\kappa_2\theta + \alpha + d\bigr] w.\eeq
Indeed, by (\ref{drate}) it follows
$$f(u) + d\,v ~=~f(u) + d u - [du - dv] ~\geq~-d (u-v).$$
{}From (\ref{wpso}) we conclude that, if $w(0,x)\geq 0$ for all $x\in\R$, then also $w(t,x)\geq 0$
for all $t\geq 0, x\in\R$.
\endproof

In  this section we focus the analysis on 

{\bf CASE 1:}  {\it The density of insects is large for $x\to +\infty$, but vanishingly small as $x\to -\infty$. 
All trees and insects are healthy in the limit as $x\to -\infty$, while they are increasingly infected as $x \to +\infty$.  }

 We seek traveling wave solutions of  (\ref{71}), having the form
\bel{TW6}u(t,x)\,=\, U(x-ct),\qquad v(t,x)\,=\, V(x-ct),\qquad \theta(t,x)\,=\,\Theta(x-ct), \qquad \alpha=\alpha(x-ct).\eeq
This leads to the system
\bel{ode6}
\left\{ \bega{rl} U'' +c U'+ f(U) - \alpha(x)\, U&=~0\,,\\[2mm]
V'' + c V'+ \kappa_2 (U-V) \Theta  - d \,V - \alpha(x)\, V&=~0,\\[2mm]
c\,\Theta'+ \kappa_1  V (1-\Theta)&=~0 \,,\enda\right.\eeq
with asymptotic conditions
\bel{asco0}
\left\{ \bega{rl} U(-\infty)&=~0\,,\\[1mm]
V(-\infty)&=~0,\\[1mm]
\Theta(-\infty)&=~0 \,.\enda\right.
\qquad\qquad \left\{ \bega{rl} U(+\infty)&=~1\,,\\[1mm]
V(+\infty)&=~V^*,\\[1mm]
\Theta(+\infty)&=~1 \,.\enda\right.
\eeq
Here $V^*=\kappa_2/( \kappa_2+d)$.

Assuming that the function $f$ satisfies {\bf (A1)},
there exists a unique speed $c^*<0$ such that the uncontrolled scalar equation
$$u_t~=~u_{xx} + f(u)$$
admits a traveling wave solution with speed $c^*$.    Moreover,  by the analysis in \cite{BCS1},
for every $c>c^*$,
there exists a non-negative control function $\alpha(\cdot)$ with minimum $\L^1$ norm,
such that the first equation in (\ref{71}) admits a traveling profile with speed $c$.

\begin{definition}\label{d:62} Let $\D$ be the domain at (\ref{PID}). 
Given an integrable function $\alpha\in \L^1(\R)$
and a constant $c<0$, we say that
the triple of functions $(U,V,\Theta): \R\mapsto \D$ is a supersolution (respectively, a subsolution) of the system (\ref{ode6}) if
\begi
\item[(i)] The functions $U,V$ are in $W^{2,1}(\R)$, i.e., they have integrable second derivatives.
\item[(ii)] The function $\Theta$ is   absolutely continuous.
\item[(iii)] The left hand sides of (\ref{ode6}) are $\leq 0$ 
(respectively: $\geq 0$) at a.e.~point $x\in \R$.
\endi
\end{definition}

Starting with a solution to the first equation, constructing a supersolution to the whole system
(\ref{ode6}) is an easy matter.

\begin{lemma}\label{l:62}
Let $u=U(x)$ be a stationary solution for the first equation in (\ref{ode6}), for some 
control $\alpha\in \L^1(\R)$ and some speed $c<0$. Define
$$\ov\Theta(x)~=~1 - \exp\left\{ {\kappa_1\over c}\int_{-\infty}^x U(y)\,dy\right\}.$$ 
Then 
the triple of functions
\bel{usol}(u^+,v^+,\theta^+)(t,x)~\doteq~\Big(U(x),\,\min\bigl\{ U(x), V^*\bigr\},\, \overline{\Theta}(x)\Big)\eeq
provides an upper solution to the system (\ref{ode6}). 
\end{lemma}

{\bf Proof.} We need to show that, by inserting the functions $u^+, v^+, \theta^+$ in (\ref{ode6}),
the left hand sides are all $\leq 0$.  A direct computation yields
$$c \theta_x^+ +\kappa_1(1-\theta^+) v^+~=~\kappa_1 (1-\ov\Theta) (v^+-U)~\leq ~0.$$
Moreover, at points where $v^+(x)=U(x)$, by (\ref{drate}),  one has
$$(v^+)'' + c (v^+)' +\kappa_2(U-v^+)\theta^+ -\alpha(x) v^+ - dv~=~U'' + cU' - \alpha(x) U - dU
~=~-f(U) - d U ~\leq~0.$$
Finally, at points where $U(x)\geq V^*$ and hence $v^+(x)=V^*$ one has
$$\bega{l}(v^+)'' + c (v^+)' +\kappa_2(U-v^+)\theta^+ -\alpha(x) v^+ - dv^+~=~
\kappa_2(U-V^*)\ov\Theta -\alpha(x) V^* - d V^*\\[2mm]
\qquad\qquad \leq~\kappa_2(1-V^*) - dV^*~=~0,\enda$$
completing the proof.\endproof

\iffalse
{\color{blue}
\begin{lemma}\label{l:n62}
Let $u=U(x)$ be a stationary solution for the first equation in (\ref{ode6}), for some 
control $\alpha\in \L^1(\R)$ and some speed $c<0$. Define
$$\ov\Theta(x)~=~1 - \exp\left\{ {\kappa_1\over c}\int_{-\infty}^x U(y)\,dy\right\}.$$ 
Then the triple of functions
\bel{lsol}
\eeq
provides an upper solution to the system (\ref{ode6}). 
\end{lemma}

{\bf Proof. }
}
\fi

In the remainder of this section we will show that the same control $\alpha(\cdot)$ yields 
a traveling profile for the system (\ref{71}), i.e.~a stationary solution to (\ref{ode6}) with 
asymptotic conditions (\ref{asco0}) as $x\to \pm\infty$. In view of Lemma~\ref{l:62},
relying on the monotonicity property of the system (\ref{ode6}), to prove the result
it remains to construct a subsolution $(u^-, v^-,\theta^-)$,
with the same asymptotic conditions (\ref{asco0}).

Introducing the variable $W=V'$, the last two equations in (\ref{ode6}) are equivalent to 
the first order system
\bel{ode8}
\left\{ \bega{rl} V'&=~W,\\[2mm]
W' &=~-c W - \kappa_2 (U-V) \Theta   +\alpha V + d \,V ,\\[2mm]
\Theta'&=\ds ~ -{\kappa_1\over c}\,  V (1-\Theta) \,,\enda\right.\eeq

Linearizing (\ref{ode8}) at the point $(U,V,W,\Theta)=(1,0,0,0)$ we obtain
 \bel{lin8}
 \begin{pmatrix} V'\cr W'\cr \Theta' \end{pmatrix} ~=~
 \begin{pmatrix} 0&1&0\cr d & -c & -\kappa_2\cr -{\kappa_1/c}&0&0 \end{pmatrix}
\begin{pmatrix} V\cr W\cr \Theta \end{pmatrix} + G(V,\Theta) + H(V,\Theta,x),\eeq
where
\bel{GH}G(V,\Theta)~=~\begin{pmatrix} 0\cr \kappa_2 V\Theta \cr \kappa_1 V\Theta /c\end{pmatrix},
\qquad\qquad   H(V,\Theta,x)~=~\begin{pmatrix} 0\cr \kappa_2\Theta \bigl(1-U(x)\bigr)
 + V\alpha(x) \cr 0 \end{pmatrix}.\eeq
The eigenvalues of the $3\times 3$  matrix in (\ref{lin8}) are the roots of the
 characteristic polynomial 
\bel{pla}
p(\lambda) ~=~ \det \begin{pmatrix}
\lambda && -1 && 0\\[1mm]
-d && c+\lambda && \kappa_2\\[1mm]
{\kappa_1 / c} && 0 && \lambda \end{pmatrix} ~=~
\lambda^3+c\lambda^2-d\lambda-{\kappa_1\kappa_2 \over c}.\eeq

\begin{figure}[ht]
\centerline{\hbox{\includegraphics[width=10cm]{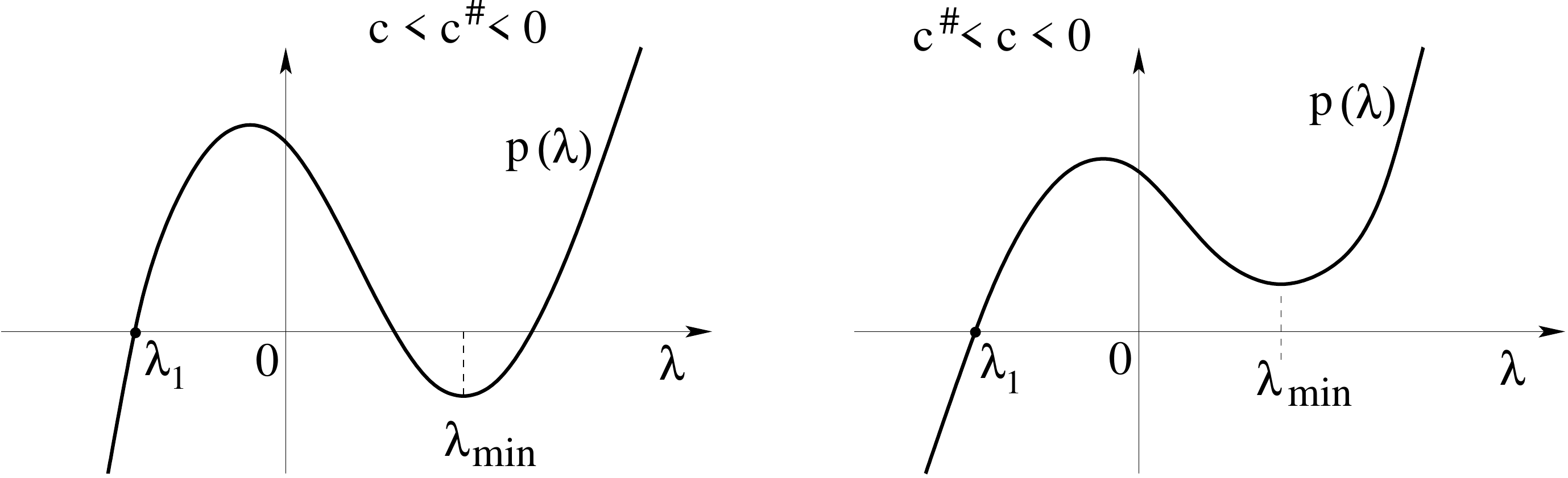}}}
\caption{\small Left: the characteristic polynomial  (\ref{pla}) in the case $c<c^\sharp$.
Right: the case $c^\sharp<c<0$.
}
\label{f:tw3}
\end{figure}
Since
$$p(0)~=~-{\kappa_1\kappa_2\over c}~>~0,$$
as shown in Fig.~\ref{f:tw3}
the polynomial $p(\lambda)$ will have  two  positive real roots 
if and only if $p(\lambda_{\min})\leq 0$, where $$\lambda_{\min}~=~{-c+\sqrt{c^2+3d}\over 3}$$ is the positive zero of $p'(\lambda)$. That is
\bel{lm}p(\lambda_{\min})~=~-{\kappa_1\kappa_2 \over c}+{cd\over 3}+{2\over 27}c^3 
- {2\over 27}(c^2+3d)^{{3\over 2}}~\leq ~0.\eeq
Differentiating the left hand side of \eqref{lm} w.r.t. $c$, we obtain
$$ {\kappa_1\kappa_2 \over c^2}~+~{d\over 3}~+~{2\over 9}c^2
~- ~{2\over 9}(c^2+3d)^{{1\over 2}}c ~ > ~ 0 \qquad \text{ for all} \quad c~<~0.
$$
Therefore, if $c^\sharp<0$ is a value for which (\ref{lm}) is satisfied as an equality,
than any value $c\leq c^\sharp$ will satisfy the inequality \eqref{lm}. 

To explicitly determine the value $c=c^\sharp$  for which the expression in (\ref{lm}) vanishes,
%\bel{lme}-{\kappa_1\kappa_2 \over c}+{cd\over 3}+{2\over 27}c^3 
%- {2\over 27}(c^2+3d)^{{3\over 2}}~= ~0.\eeq
we move the last term to the right side, square both sides and simplify the equation to get
$$ {c^2 d^2 \over 27} + {4 d^3 \over 27} + {4\over 27} \kappa_1\, \kappa_2\, c^2 + {2\, \kappa_1\, \kappa_2\, d\over 3} - {\kappa_1^2 \kappa_2^2 \over c^2}~=~0.
$$
We solve the above equation for $c^2$, and take the negative square root.
This yields
\bel{cstar}
c^\sharp~\doteq~-\bigg({-2 d^3-9\kappa_1\kappa_2 d + 2\big(d^2 + 3\kappa_1\kappa_2\big)^{3/2}\over d^2+4 \kappa_1\kappa_2}\bigg)^{1/2}.\eeq
{}From the above analysis it follows
\begin{lemma}\label{l:71} For $c^\sharp<c<0$, the $3\times 3$ Jacobian matrix at (\ref{lin8}) has one
negative eigenvalue and two complex conjugate eigenvalues, with positive real part.
\end{lemma}
Calling 
\bel{eigi} \lambda_1\,<\,0, \qquad \lambda_2\,=\, a+ib,\qquad \lambda_3\,=\, a-ib,
\eeq the three eigenvalues, with $a,b>0$, we obtain
three
corresponding eigenvectors:
\bel{evi} \bfv_i~ =~ \begin{pmatrix}
1  \\[1mm]
\lambda_i \\[1mm]
-{ \kappa_1 /c \lambda_i } \end{pmatrix},\qquad\qquad i=1,2,3. \eeq
Notice that $\bfv_1$ has real entries, while $\bfv_2,\bfv_3$ are complex valued.
Taking the real and imaginary parts, we obtain the two vectors 
\bel{2wv}\bfw_2~=~\begin{pmatrix} 1\cr a\\[1mm]\ds -{\kappa_1 a\over c(a^2 + b^2)}\end{pmatrix},\qquad\qquad 
\bfw_3~=~\begin{pmatrix} 0\cr b\\[1mm] \ds{\kappa_1 b\over c(a^2 + b^2)}\end{pmatrix},\eeq
which satisfy
\bel{Sigma} \Sigma~\doteq~\hbox{span}\{\bfw_2, \bfw_3\} ~=~ \hbox{span}\{ \bfv_2,\bfv_3\}.\eeq
In particular, a direct computation shows that the linear system
 \bel{lin7}
 \begin{pmatrix} V'\cr W'\cr \Theta' \end{pmatrix} ~=~
 \begin{pmatrix} 0&1&0\cr d & -c & -\kappa_2\cr -{\kappa_1/c}&0&0 \end{pmatrix}
\begin{pmatrix} V\cr W\cr \Theta \end{pmatrix} \eeq
admits the solution
\bel{sol1}\begin{pmatrix} V\cr W\cr \Theta \end{pmatrix}(x) ~=~ 
Ae^{ax}\begin{pmatrix} \sin bx\cr b\cos bx+a\sin bx \cr {\kappa_1 \over c (a^2+b^2)}\bigl(b\cos bx - a\sin bx\bigr) \end{pmatrix}.
\eeq
We can now prove the main result of this section, on the existence of controlled traveling waves
for Model~2.

\begin{theorem}\label{t:61} Let $f$ satisfy the assumptions {\bf (A1)} together with (\ref{drate}).
Let $c^\sharp<c<0$ and let $U:\R\mapsto [0,1]$ be an increasing solution to the first equation in (\ref{ode6}),
with asymptotic conditions as in (\ref{asco0}), for some nonnegative control function $\alpha\in \L^1(\R)$ with bounded support.
Then there exist solutions $V, \Theta$ of the remaining two equations in (\ref{ode6}), 
with asymptotic conditions (\ref{asco0}).
\end{theorem}

\begin{figure}[ht]
\centerline{\hbox{\includegraphics[width=9cm]{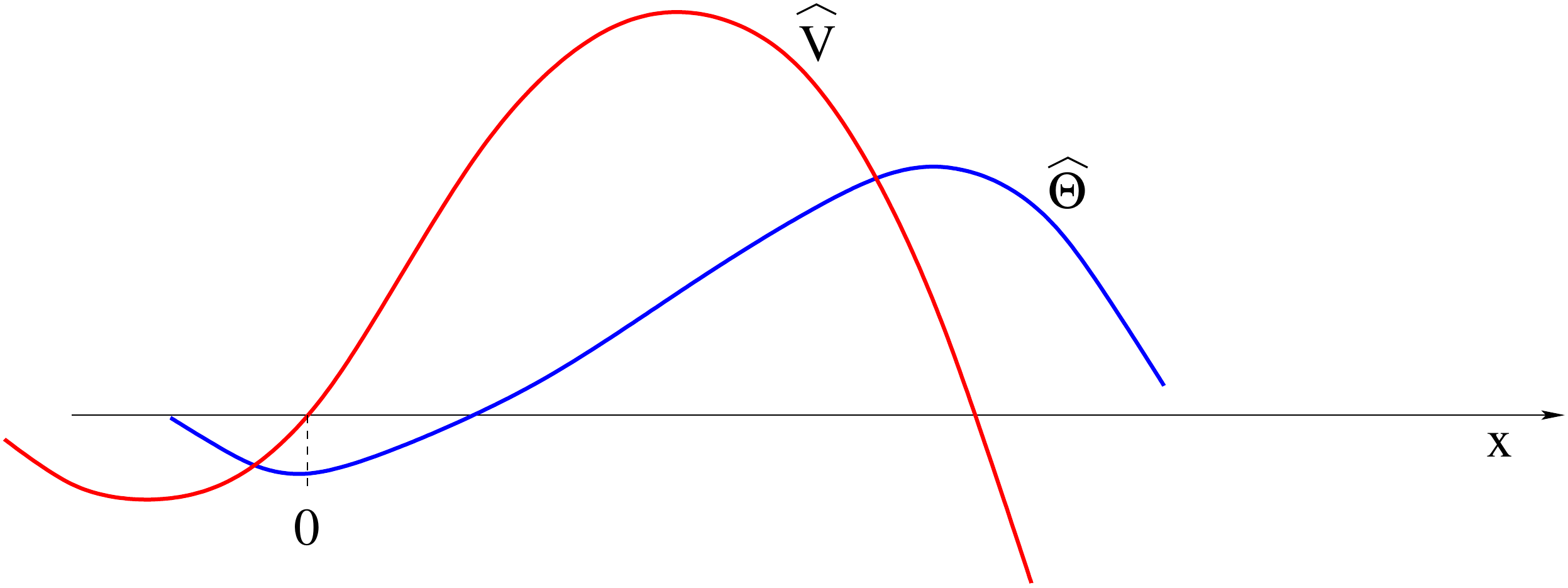}}}
\caption{\small  A particular solution (\ref{VT}) to the linear system (\ref{lin7}).}
\label{f:co47}
\end{figure}

{\bf Proof.}  By Lemma~\ref{l:62} we already have an upper solution of (\ref{ode6}) satisfying the 
asymptotic conditions (\ref{asco0}). It remains to construct a lower solution.
\v
{\bf 1.} Let
 $\varphi_0\in ]0, \pi/2[$ be the angle such that
$$\cos\varphi_0 ~=~{ a \over  \sqrt{a^2+b^2}}\,,\qquad  \sin\varphi_0~={ b \over  \sqrt{a^2+b^2}}\,.$$
Then by (\ref{sol1}) the functions
\bel{VT}
\Hat{V}(x)~=~e^{ax} \sin bx, \qquad\qquad \Hat{\Theta}(x)~=~- {\kappa_1\over c\sqrt{a^2+b^2}} e^{ax}
\sin \big[b x - \vp_0\big]\eeq
%\bel{VT}
%\widecheck{V}(x)~=~A e^{ax} \sin b(x-\theta_0), \qquad\qquad \widecheck{\Theta}(x)~=~- {\kappa_1\over c\sqrt{a^2+b^2}} Ae^{ax}
%\sin \big[b(x-\theta_0) - \vp_0\big],\eeq
provide one particular  solution to the linear system (\ref{lin7}), as shown in Fig.~\ref{f:co47}.
%, for any choice of  $A, \theta_0$.
\v
{\bf 2.} 
Given $x_0\in \R$, for any $\ve>0$, call $(V_\ve, W_\ve, \Theta_\ve)$ the solution to the system
(\ref{ode6}) with initial data
\bel{ida5}
V_\ve(x_0)\,=\,0, \qquad W_\ve(x_0)\,=\,V_\ve'(x_0)\,=\, \ve \Hat V'(0),\qquad
\Theta_\ve(x_0)\,=\,\ve \Hat \Theta(0),\eeq
in the special case where $\alpha(x)=0$ and $U(x)=1$ for all $x$.
By standard ODE theory, as $\ve\to 0$ we have the  convergence
\bel{ufc}
\ve^{-1} V_\ve(x+x_0)~\to~\Hat V(x),\qquad \ve^{-1} V'_\ve(x+x_0)~\to~\Hat V'(x),
\qquad \ve^{-1} \Theta_\ve(x+x_0)~\to~\Hat \Theta(x),\eeq
uniformly for $x$ in bounded intervals.

Since we are assuming that the control $\alpha(\cdot)$ has bounded support, for any $\epsilon_0>0$
we can choose  $x_0>0$ large enough so that
\bel{esmall} \alpha(x)\,=\,0,\qquad\qquad
1-\epsilon_0\,\leq U(x)\,\leq  1,\qquad\forall x\geq x_0\,.\eeq
By choosing $\ve,\epsilon_0>0$ small enough, we obtain an exact solution
of (\ref{ode6}) on an interval $[x_0, x_1]$, with  $x_1-  x_0\leq 4\pi/b$ and $\alpha, U$ as in (\ref{esmall}), such that
\bel{VTP}
\left\{ \bega{l} V_{\ve}(x_0)=0,\cr
V_{\ve}(x_1)=0,\enda\right.\qquad\quad  \left\{ \bega{l} \Theta_{\ve}(x_0) <0,\cr
\Theta_{\ve}(x_1)>0,\enda\right.\qquad \quad \left\{\bega{l}
V_{\ve}(x)> 0\quad \hbox{for}~~ x_0<x<x_1\,,\cr
V_{\ve}'(x_1)<0\,.\enda\right.\eeq

Restricted to the half line $\,]-\infty, x_1]$, our lower solution is then defined as
\bel{ls1}
v^-(x)~=~\left\{ \bega{cl} 0\quad &\hbox{if} ~~x<x_0\,,\cr
V_{\ve}(x)\quad &\hbox{if} ~~x\in [x_0, x_1]\,,\enda\right.
\qquad\quad 
\theta^-(x)~=~\left\{ \bega{cl} 0\quad& \hbox{if} ~~x<x_0\,,\cr
\max\bigl\{ \Theta_{\ve}(x),0\bigr\} \quad &\hbox{if} ~~x\in [x_0, x_1]\,.\enda\right.\eeq

\v
{\bf 3.} Next, we extend this subsolution to the remaining half line $[x_1,\,+\infty[\,$. 

As a first step, we define the constant function
$$\Tilde \theta(x)~\doteq ~\Theta_\ve(x_1)~>~0,$$
and let $\Tilde v$ be the solution to 
\bel{v3p}v''~=~-c v' -\kappa_2(1-\epsilon_0- v) \tilde \theta_1  + dv,\eeq
on the domain $x\in [x_1, +\infty[\,$, with boundary conditions
\bel{v4p}v(x_1)\,=\,0, \qquad\qquad v(+\infty) \,=\, V^\dagger\,\doteq\,{\kappa_2(1-\epsilon_0) \Theta_{\ve}(x_1) \over \kappa_2\Theta_{\ve}(x_1) + d}\,.\eeq
An explicit computation yields
$$\Tilde v(x) ~=~ V^{\dagger} \big(1 - e^{\lambda_0(x-x_1)}\big), \qquad\quad  \lambda_0 \,=\,
 {-c-\sqrt{c^2+4(\kappa_2\Theta_{\ve}(x_1)+d)} \over 2}.$$
 Notice that, for $x>x_1$, the couple $(\Tilde v, \Tilde\theta) $ provides a lower solution
 to the last two equations in (\ref{ode6}).   However, this subsolution does not yet satisfy
 the asymptotic conditions in (\ref{asco0}).  One more step is thus needed.
\v
{\bf 4.}
 For  $x>x_1$ we let $\theta^-$ be the solution to
\bel{th'}\theta'~=~{-\kappa_1\over c} \Tilde v(x) (1-\theta),\qquad\qquad \theta(x_1) = \Theta_{\ve}(x_1),\eeq
where $\Tilde v$ is the function constructed in the previous step.  More explicitly, this means
$$\theta^-(x)~=~1-\big(1-\Theta_\ve(x_1)\big)\exp\left\{\int_{x_1}^x {\kappa_1\over c}\Tilde v(z)\, dz
\right\}.
$$
Observe that, since $c<0$ and $\Tilde v_1(x)\to V^\dagger >0$ as $x\to +\infty$, the above solution $\theta^-$ is monotone increasing and satisfies $\theta^-(x)\geq \Tilde \Theta_\ve(x_1)$ as $x \in [x_1, +\infty)$ and $\theta^-(x)\to 1$ as $x\to +\infty$.

We then define  $v^-$ to be the solution of
\bel{v5p}v''~=~-c v' -\kappa_2(U(x)- v) \theta^-(x)  + dv,\eeq
on the domain $x\in [x_1, +\infty[\,$, with boundary conditions
\bel{v6p}v(x_1)\,=\,0, \qquad\qquad v(+\infty) \,=\, V^*\,\doteq\,{\kappa_2\over\kappa_2+d}\,.\eeq
Observing that 
$$U(x)\,\geq\, 1-\epsilon_0\,,\qquad\qquad \theta^-(x)\,\geq\,\Tilde\theta\qquad\qquad\forall x\geq x_1\,,$$
by a comparison argument we conclude
$$v^-(x)~\geq~\Tilde v(x)\qquad\qquad \forall x\geq x_1\,.$$
It is now clear that the couple $(v^-, \theta^-)$ provides a subsolution, restricted 
to the half line $[x_1, +\infty[\,$.
Since $v^-(x)\geq 0$ for all $x\in\R$ while $v^-(x_1)=0$, it follows that at the 
junction point $x_1$ the left and right derivatives of $v^-$ satisfy
\bel{lrdv}(v^-)'(x_1-)~\leq~0~\leq~(v^-)'(x_1+).\eeq
Hence $(v^-, \theta^-)$ is a subsolution defined on the whole real line, which satisfies all the
asymptotic conditions in (\ref{asco0}).

\begin{figure}[ht]
\centerline{\hbox{\includegraphics[width=13cm]{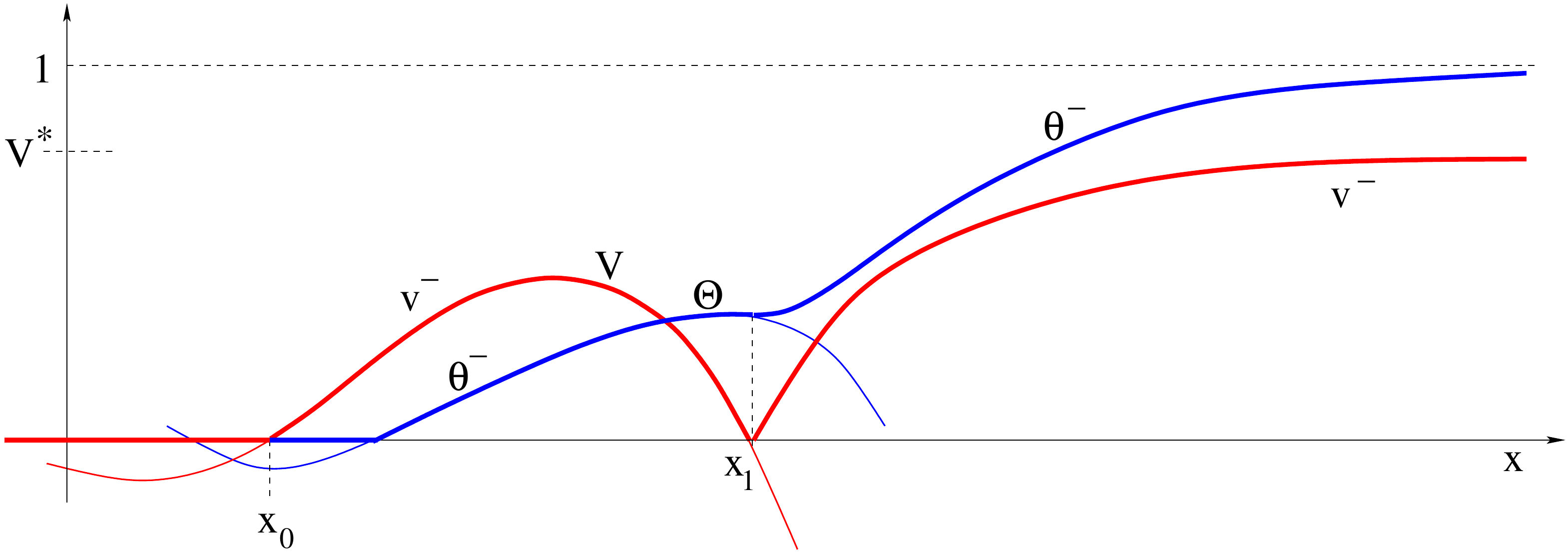}}}
\caption{\small The lower solution constructed  at (\ref{ls1}) and at (\ref{th'})-(\ref{v6p}), separately on the half lines where $x\leq x_1$ and $x\geq x_1$.}
\label{f:co46}
\end{figure}

\v
{\bf 5.} Having constructed a supersolution and a subsolution of (\ref{ode6}) with
\bel{SSS}v^-(x)\leq v^+(x),\qquad \qquad \theta^-(x)\leq \theta^+(x)\qquad\forall x\in \R,\eeq
the existence of an exact solution follows by a standard monotonicity argument.
Namely, since the (component-wise) supremum 
of two subsolutions is also a subsolution, we can define
$$\bigl(V(x), \Theta(x)\bigr)~=~\sup_{(v^-, \theta^-)\in {\cal S}} ~ \bigl(v^-(x), \theta^-(x)\bigr),$$
where the supremum is taken over the set ${\cal S}$ of all subsolutions which satisfy (\ref{SSS}).
More precisely:
$$v^-(x)~\leq~\min\bigl\{ U(x), V^*\bigr\},\qquad\qquad \theta^-(x)~\leq~\ov\Theta(x).$$
By construction, our subsolutions and supersolutions all satisfy the same asymptotic conditions
at (\ref{asco0}).   Hence the same holds for the exact solution.
\endproof

\section{Nonexistence of controlled traveling profiles with slow speed}  
\label{sec:7}
\setcounter{equation}{0}
In this section we continue the analysis of the system (\ref{71}), focusing on

{\bf CASE 2:} The density of insects is large  for $x\to +\infty$  as well as for $x\to -\infty$. 
Insects and trees
are all healthy in the limit as $x\to -\infty$, while they are increasingly infected as $x \to +\infty$.

We consider the possibility of  using a control $\alpha(\cdot)$
to reduce the density of insects in the intermediate region  between the  healthy and contaminated zone.
In principle, this should provide a ``buffer zone", separating the healthy population from the sick one, 
thus slowing down the spread of the contamination.
Our analysis, however, will show that this strategy is not effective.  Namely, it cannot yield any
traveling wave profile with slower propagation speed.

To state a precise result in this direction,
we first study the asymptotic behavior of a traveling wave
as $x\to -\infty$.
To fix ideas, let a speed $c<0$ be given.
We seek a control $\alpha\in \L^1(\R)$ and a solution of (\ref{71}) in the form of a traveling wave
(\ref{TW6}). This leads again to the system (\ref{ode6}). However, the
asymptotic conditions (\ref{asco0}) are now replaced by
\bel{asco}
\left\{ \bega{rl} U(-\infty)&=~1\,,\\[1mm]
V(-\infty)&=~0,\\[1mm]
\Theta(-\infty)&=~0 \,,\enda\right.
\qquad\qquad \left\{ \bega{rl} U(+\infty)&=~1\,,\\[1mm]
V(+\infty)&=~V^*,\\[1mm]
\Theta(+\infty)&=~1 \,.\enda\right.
\eeq

\begin{theorem}\label{t:71}  Let $c^\sharp$ be the constant in (\ref{cstar}), and consider any speed $c$ with $c^\sharp<c<0$.
Then the system (\ref{ode6}) does not admit any solution 
$x\mapsto \bigl(U(x),V(x),\Theta(x)\bigr)\in\D$ with asymptotic conditions (\ref{asco}), 
for any control
$\alpha\in \L^1(\R)$.
\end{theorem}

{\bf Proof.}
The proof will be achieved by showing that, even by adding a control $\alpha\in\L^1(\R)$
in the equations (\ref{ode8}), one cannot achieve solutions  such that $V(x), \Theta(x)$
converge to zero
as $x\to -\infty$, and satisfy the constraint $V(x),\Theta(x)\in [0,1]$ for all $x\in\R$.
The argument will be given in several steps.

\v
{\bf 1.} Assume that, on the contrary, a traveling wave solution $(U,V,\Theta)$, exists, with the 
prescribed asymptotic behavior as $x\to -\infty$.   A contradiction will be obtained by showing that 
the control $\alpha(\cdot)$ cannot be integrable. 

As a preliminary, we observe that the assumption $\alpha\in \L^1(R)$ implies
that the traveling wave profile $U(\cdot)$, i.e. the solution to 
\bel{TVU}U'' + cU' +f(U) - \alpha(x) U~=~0,\qquad\qquad U(-\infty) ~=~U(+\infty)~=~1,\eeq
 satisfies
\bel{UL1}\int_{-\infty}^0 \bigl(1-U(x)\bigr)\, dx ~<~+\infty. \eeq

On the space $\R^3$, it will be convenient to use a new system of coordinates
$y=(y_1,y_2,y_3)$ corresponding to the basis $\{ \bfv_1, \bfw_2,\bfw_3\}$ defined in \eqref{evi}, \eqref{2wv}.
Let $x\mapsto Y(x)= \bigl( Y_1(x), Y_2(x), Y_3(x)\bigr)$ be the coordinates of the traveling profile
$(V,W,\Theta)$ w.r.t.~this new basis. By construction, the system (\ref{ode8})
can be written as
\bel{ode9} \begin{pmatrix} Y_1'\cr Y_2'\cr Y_3' \end{pmatrix} ~=~
 \begin{pmatrix} \lambda_1&0&0\cr 0 & a & -b\cr 0&b&a \end{pmatrix}
\begin{pmatrix} Y_1\cr Y_2\cr Y_3 \end{pmatrix} + \Tilde G(Y)+\Tilde H(Y,x) ,\eeq
where, in view of (\ref{GH}) and (\ref{UL1}),
 the nonlinear perturbations $\Tilde G, \Tilde H$ satisfy the bounds
\bel{ghb}\bigl| \Tilde G(Y)\bigr|~\leq~C_0 \,|Y|^2,\qquad\qquad \bigl| \Tilde H(Y,x)\bigr|~\leq~C_0 \,|Y|\, \beta(x),\eeq
for some constant $C_0$ and some integrable function $\beta\in \L^1\bigl(]-\infty, 0]\bigr)$.
\v
{\bf 2.} By (\ref{ode9})-(\ref{ghb}) there exists a constant $C_1$ such that
$$\left|{d\over dx} \bigl|Y(x)\bigr|~\right|~\leq~C_1\Big( 1+\bigl|\beta(x)\bigr|\Big)  \bigl|Y(x)\bigr|.$$
Since $\beta\in \L^1$, we conclude that the vector $Y(x)$ cannot vanish at any point
$-\infty<x\leq 0$.
\v
{\bf 3.} Introducing the radius 
$r(x)=  \bigl| Y(x)\bigr|$, we now consider the normalized vector $\xi$, such that
$$\xi(x)~=~(\xi_1, \xi_2,\xi_3)(x)~=~{Y(x)\over \bigl| Y(x)\bigr|}\,,\qquad\qquad Y (x)\,=\, r(x)\xi(x)\,.$$
By (\ref{ode9}), denoted the $3\times 3$ matrix in \eqref{ode9} as $A$, this vector $\xi$ satisfies
\bel{li2}
\xi'(x)~=~A \xi + \tilde g(r,\xi) + \tilde h(r,\xi,x) - \Big\langle A \xi + \tilde g(r,\xi) + \tilde h(r,\xi,x) ~,~\xi
\Big\rangle \xi\,,\eeq
where
\bel{tgh}\tilde g(r,\xi)~=~r^{-1} \Tilde G(r\xi),\qquad\qquad \tilde h(r,\xi,x)~=~r^{-1} \Tilde H(r\xi,x).\eeq
Since  $r(x)\to 0 $ as $x\to -\infty$,  by (\ref{ghb}) and 
(\ref{tgh}) we have
$$\bigl|\tilde g(r,\xi)\bigr| \leq C_0 \bigl|r(x)\bigr|,
\qquad  \bigl|h(r,\xi,x)\bigr|~\leq~C_0\bigl|\beta(x)\bigr|\qquad\qquad \lim_{x\to -\infty}~ \Big|\tilde g\bigl(
r(x),\xi\bigr)\Big|~=~0 ,$$
uniformly for all $|\xi|=1$.
\v
{\bf 4.} Based on the previous step, we observe that, as $x\to -\infty$, the evolution of the normalized vector  $\xi=(\xi_1,\xi_2,\xi_3)$ satisfies an equation of the form
\bel{li3}
\xi'~=~A \xi - \bigl\langle A \xi\,,~\xi
\bigr\rangle \xi + g(x) + h(x),\eeq
where $h\in \L^1$ while $\lim_{x\to -\infty}g(x)=0$.
We claim that, as $x\to -\infty$, %and $r(x) = \bigl|Y(x)\bigr|\to 0$, 
two cases are possible

{\bf Case 1:} $\xi_1(x)\to \pm 1$.

{\bf Case 2:} $\xi_1(x)\to 0$.

 Indeed, by (\ref{li2}),  the first component of  the vector $\xi$ satisfies the ODE
\bel{xi1}\xi_1'(x) ~=~ (\lambda_1-a)\xi_1(1-\xi_1^2) +g_1(x) + h_1(x),  
\eeq
where $h_1\in \L^1$ while $\lim_{x\to -\infty}g_1(x)=0$.

For any $\delta\in \,]0, 1/2]$, consider the set
\bel{Iddef}\bega{rl} 
I_\delta&=~\ds \bigg\{ \bar x\leq 0\,;\quad \bigl| \xi_1(\bar x)\bigr|\geq\delta,\quad 
\int_{-\infty}^{\bar x} \bigl|h_1(y)\bigr|\, dy< {\delta\over 2}\,,
\\[4mm]
&\ds\qquad\qquad \qquad 
\bigl|g_1(x)\bigr| < (a-\lambda_1) {\delta(1-\delta^2)\over 2}  \qquad\forall x\leq \bar x\bigg\}.\enda\eeq
Assume that one of these sets $I_\delta$ is nonempty, say $\bar x\in I_\delta$.
We claim that 
\bel{x1b}\bigl|\xi_1(x)\bigr|~>~{\delta\over 2}\qquad\forall x\leq \bar x.\eeq
Indeed, consider the function
\bel{phdef}\phi(x)~\doteq~\bigl|\xi_1(x) \bigr|-  \int_{-\infty}^ x \bigl|h_1(y)\bigr|\, dy.\eeq
Recalling that $\lambda_1<0<a$, for $x\leq \bar x$ by (\ref{Iddef}) we have the implication
$$|\xi_1|\in \left[{\delta\over 2}, \delta\right]\quad\implies\quad
\phi'(x)~\leq ~(\lambda_1-a)|\xi_1|(1-\xi_1^2) +\bigl|g_1(x)\bigr| ~<~0.$$
If there exist $x_1<x_2\leq\bar x$ such that
\bel{xi123}{\delta\over 2} ~=~\bigl|\xi_1(x_1)\bigr| ~<~\bigl|\xi_1(x_2)\bigr|~=~\delta,\eeq 
then 
$${\delta\over 2} ~=~\bigl|\xi_1(x_2)\bigr|-\bigl|\xi_1(x_1)\bigr|~=~\phi(x_2) - \phi(x_1)
+\int_{x_1}^{x_2}
\bigl|h_1(y)\bigr|\, dy~<~\int_{-\infty}^{\bar x}\bigl|h_1(y)\bigr|\, dy~<~{\delta\over 2}\,,$$
reaching a contradiction.

Using (\ref{x1b}), we now show that
\bel{limp}
\lim_{x\to -\infty} \bigl|\xi_1(x)\bigr|~=~\lim_{x\to -\infty} \phi(x)~=~1.\eeq
Indeed, the first identity is an immediate consequence of (\ref{phdef}).
To prove the second equality, let any $\ve\in \,]0, \delta/4]$ be given.
Choose $x^*< \bar x$ such that 
$$\int_{-\infty}^{x^*} \bigl|h_1(y)\bigr|\, dy< \ve,\qquad \qquad 
\bigl|g_1(x)\bigr| < \ve  \qquad\forall x\leq x^*.$$
Observing that
$${\delta\over 4}~\leq~\bigl|\xi_1(x)\bigr|-\ve~\leq~\phi(x)~\leq~|\xi_1(x)\bigr|\qquad\forall x\leq x^*,$$
{}from (\ref{xi1}) we obtain
$$\phi'(x)~\leq~(\lambda_1-a) {\delta\over 4} \bigl( 1 - \phi^2(x)\bigr) + \ve ~<~0,$$
where the last inequality holds as long as 
$$1 - \phi^2(x)~\geq~{4\ve\over (a-\lambda_1)\delta}\,.$$
We thus conclude
$$\liminf_{x\to -\infty} \bigl(1 - \phi^2(x)\bigr)~\leq~{4\ve\over (a-\lambda_1)\delta}\,.$$
Since $\ve>0$ can be arbitrarily small, this yields the second identity in (\ref{limp}).
\v
The previous analysis has shown that, if one of the sets $I_\delta$ is nonempty, then 
(\ref{limp}) holds, hence Case 1 occurs.

The remaining possibility is that all sets $I_\delta$ are empty. In this case, for every 
$\delta>0$ we can find $x^*<0$ such that 
$$
\int_{-\infty}^{x^*} \bigl|h_1(y)\bigr|\, dy< {\delta\over 2}\,,
\qquad\qquad
\bigl|g_1(x)\bigr| < (a-\lambda_1) {\delta(1-\delta^2)\over 2}  \qquad\forall x\leq x^*.$$
This implies $\bigl|\xi_1(\bar x)\bigr|< \delta$ for all $\bar x\leq x^*$, otherwise $\bar x\in I_\delta$
against the assumption. We thus conclude that Case 2 holds true.
\v
{\bf 5.} We show that Case 1 leads to a contradiction.
Indeed, given $\ve>0$, by choosing $x_0<\!<0$ we achieve
\bel{Y1big}|Y(x)|\leq\ve,\qquad \int_{-\infty}^x \bigl|\beta(x)\bigr|\, dx~\leq~\ve, 
\qquad |Y(x)|~\leq ~2 |Y_1(x)|,\qquad\forall x<x_0\,.\eeq
In this case, for any $x_1<x<x_0$ we have
$$\bigl|Y_1(x)\bigr|~\leq~e^{\lambda_1(x-x_1)}\bigl|Y_1(x_1)\bigr| +C \int_{x_1}^x e^{\lambda_1(x-y)} \Big(
\bigl|Y_1(y)\bigr|+\bigl|\beta(y)\bigr|\Big)
\bigl|Y_1(y)\bigr|\, dy.$$
Letting $x_1\to -\infty$ we obtain
\bel{ctrac}\bega{rl} \bigl|Y_1(x)\bigr|&\ds\leq~C \int_{-\infty}^x e^{\lambda_1(x-y)} \Big(
\bigl|Y_1(y)\bigr|+\bigl|\beta(y)\bigr|\Big)
\bigl|Y_1(y)\bigr|\, dy\\[4mm]
&\ds \leq~\ds C\ve  \int_{-\infty}^x e^{\lambda_1(x-y)} \Big(
\ve+\bigl|\beta(y)\bigr|\Big)
\, dy\\[3mm]
& \ds \leq~C{\ve^2\over | \lambda_1|} e^{\lambda_1 x}  + C\ve \int_{-\infty}^x \lambda_1 e^{\lambda_1 (x-y)} \left(
\int_y^x |\beta(z)|\, dz\right) dy\\[3mm]
&\leq ~\ds C\ve \Big( C_1 \ve + C_2 \int_{-\infty}^x |\beta(z)|\, dz\Big)~\leq~{\ve\over 4}\,,
\enda\eeq
provided that $\ve>0$ is chosen sufficiently small.   By the third inequality on (\ref{Y1big})
we conclude
$\bigl|Y(x)\bigr|\leq \ve/2$ for all $x<x_0$.

Iterating this argument, we obtain
$\bigl|Y(x)\bigr|\leq 2^{-k}\ve$ for every $k\geq 1$, hence $Y(x)=0$ for all $x\in \,]-\infty, x_0]$, 
reaching a contradiction.
\v
{\bf 6.} We now show that Case 2 also leads to a contradiction.
By step {\bf 3},  the last two components 
satisfy an ODE of the form
\bel{ode23}
 \begin{pmatrix} \xi_2'\cr \xi_3' \end{pmatrix} ~=~
 \begin{pmatrix} 0 & -b\cr b&0 \end{pmatrix}
\begin{pmatrix} \xi_2\cr \xi_3 \end{pmatrix} +  \phi_1(x) + \phi_2(x),\eeq
with 
$$\lim_{x\to -\infty} \phi_1(x)~=~0,\qquad\qquad \int_{-\infty}^0 \phi_2(x)\, dx ~<~+\infty.$$

On the plane $\Sigma$ at (\ref{Sigma}), it will be convenient to use polar coordinates
$(r, \vartheta)$. More precisely, by (\ref{ode23}) the evolution of the angle variable has the form
 \bel{dvt}
 {d\over dx} \vartheta(x)~=~b + \Tilde \phi_1(x) + \Tilde \phi_2(x) .\eeq
 where 
 $$\Tilde \phi_1(x)\to 0,\qquad\qquad \int_{-\infty}^{x_0} \bigl| \Tilde \phi_2(x)\bigr|\, dx~<~+\infty.$$
 
\begin{figure}[ht]
\centerline{\hbox{\includegraphics[width=16cm]{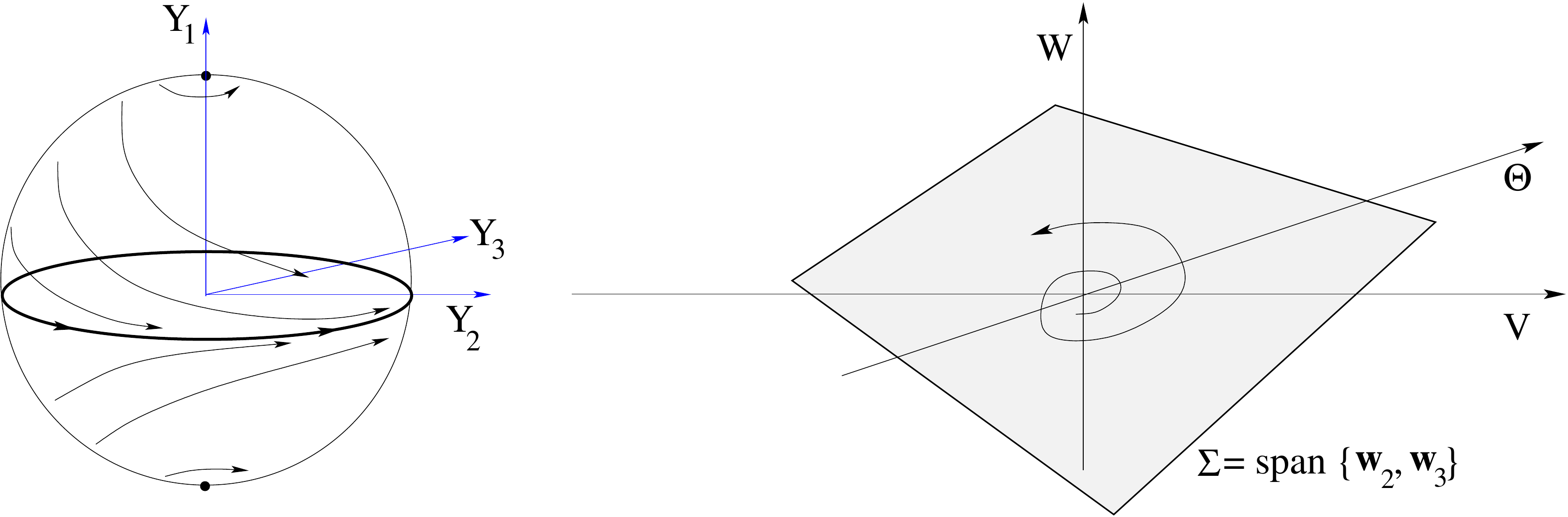}}}
\caption{\small Left: the dynamics of the unit vector $\xi= Y/|Y|$, on the surface of the unit ball
in $\R^3$.
Right: on the plane $\Sigma= \hbox{span}\{\bfw_1,\bfw_2\}$,
for certain values of the angular component $\vartheta$, 
the point $P$ with polar coordinates $(r,\vartheta)$ lies outside the admissible set
where $V\geq 0$ and $\Theta\geq 0$.
}
\label{f:tw4}
\end{figure}

 As shown in Fig.~\ref{f:tw4}, this implies that the trajectory makes infinitely many loops
around the origin, close to the plane $\Sigma$.   But this is impossible, 
because in this case, for some values of $x$, one of the components $V(x)$, $\Theta(x)$
must be negative. This concludes the proof of the theorem.
\endproof

 \end{document}